\documentclass[11pt,a4paper]{article}

 \usepackage{geometry}
 \newgeometry{ hmargin={20mm,20mm}}

\usepackage{algorithm}
\usepackage{amsmath}
\usepackage{amssymb}
\usepackage{bm}
\usepackage{booktabs}
\usepackage{caption,subcaption}
\usepackage{listings}
\usepackage{pythonhighlight}

\usepackage{natbib}

\usepackage{tikz}
\usetikzlibrary{arrows.meta,
                chains,
                positioning,
                shapes.geometric
                }
              \usetikzlibrary{shapes,positioning}

\DeclareMathOperator*{\argmin}{argmin}
\DeclareMathOperator*{\minimize}{minimize}
\DeclareMathOperator*{\maximize}{maximize}

\DeclareMathOperator{\diag}{diag}
\DeclareMathOperator{\Diag}{Diag}
\DeclareMathOperator{\Tr}{Tr}

\DeclareMathOperator{\bone}{\bf 1}
\DeclareMathOperator{\bzero}{\bf 0}
\DeclareMathOperator{\bA}{\bf A}

\DeclareMathOperator{\blb}{\bf b}

\DeclareMathOperator{\bdA}{\bf dA}
\DeclareMathOperator{\bdb}{\bf db}
\DeclareMathOperator{\bdG}{\bf dG}
\DeclareMathOperator{\bdh}{\bf dh}

\DeclareMathOperator{\bdz}{\bf dz}
\DeclareMathOperator{\bD}{\bf D}
\DeclareMathOperator{\bd}{\bf d}
\DeclareMathOperator{\bE}{\bf E}
\DeclareMathOperator{\bF}{\bf F}
\DeclareMathOperator{\bG}{\bf G}

\DeclareMathOperator{\bI}{\bf I}
\DeclareMathOperator{\bh}{\bf h}
\DeclareMathOperator{\bl}{\bf l}

\DeclareMathOperator{\bp}{\bf p}
\DeclareMathOperator{\bQ}{\bf Q}
\DeclareMathOperator{\bq}{\bf q}

\DeclareMathOperator{\br}{\bf r}

\DeclareMathOperator{\bu}{\bf u}
\DeclareMathOperator{\bV}{\bf V}
\DeclareMathOperator{\bv}{\bf v}
\DeclareMathOperator{\bW}{\bf W}
\DeclareMathOperator{\bw}{\bf w}

\DeclareMathOperator{\bx}{\bf x}
\DeclareMathOperator{\bY}{\bf Y}
\DeclareMathOperator{\by}{\bf y}
\DeclareMathOperator{\bz}{\bf z}

\DeclareMathOperator{\bbeta}{\bm \beta}

\DeclareMathOperator{\bdlambda}{ \bd\boldsymbol{\lambda} }
\DeclareMathOperator{\bdetta}{ \bd\boldsymbol{\eta}}

\DeclareMathOperator{\betta}{\bm \eta}
\DeclareMathOperator{\bgamma}{\bm \gamma}
\DeclareMathOperator{\blambda}{\bm \lambda}
\DeclareMathOperator{\bmu}{\bm \mu}
\DeclareMathOperator{\bnu}{\bm \nu}

\DeclareMathOperator{\btheta}{\bm \theta}

\providecommand{\keywords}[1]{\textbf{\textit{Keywords:}} #1}

\newcommand{\toi}[2][i]{%
  \mathop{
    \mathrm{#2}^{( #1 )}
  }
}

\newcommand{\toit}[2][i]{%
  \mathop{
    \mathrm{#2}^{{T}^{( #1 )}}
  }
}

\newcommand{\toinv}[2][i]{%
  \mathop{
    \mathrm{#2}^{{-1}^{( #1 )}}
  }
}

\begin{document}

\title{Data-driven integration of norm-penalized mean-variance portfolios}
\author{Andrew Butler and Roy H. Kwon \\ University of Toronto\\Department of Mechanical and Industrial Engineering}
\maketitle

\begin{abstract}
Mean-variance optimization (MVO) is known to be sensitive to estimation error in its inputs.  Norm penalization of MVO programs  is a regularization technique that can mitigate the adverse effects of estimation error. We augment the standard MVO program with a convex combination of parameterized $L_1$ and $L_2$-norm penalty functions. The resulting program is a parameterized quadratic program (QP) whose dual is a box-constrained QP. We make use of recent advances in neural network architecture for differentiable QPs and present a data-driven framework for optimizing parameterized norm-penalties to minimize the downstream MVO objective. We present a novel technique for computing the derivative of the optimal primal solution with respect to the parameterized $L_1$-norm penalty by implicit differentiation of the dual program. The primal solution is then recovered from the optimal dual variables. Historical simulations using US stocks and global futures data demonstrate the benefit of the data-driven optimization approach.
\end{abstract}


\maketitle

\section{Introduction}\label{sec:intro}
The \cite{Markowitz1952} mean-variance optimization (MVO) requires as input estimates of asset mean returns and covariances. Asset returns, however, are nonstationary and exhibit statistically insignificant auto-correlation \citep{Drees2002,Engle1982,Starica2005} and are therefore estimated with a high degree of uncertainty \cite{Michaud2008b}. Moreover, estimating a covariance matrix of $N$ assets requires computing $N(N+1)/2$ cross-covariance terms and therefore when the number of assets is large relative to the relevant number of historical observations then the estimated covariance matrix can be numerically unstable \citep{Ledoit2003}. As such, the sample covariance matrix and the action of its inverse is also contaminated with estimation error \citep{Jobson1980}. 

In this paper, we address the adverse effects of estimation error by augmenting the mean-variance objective with a norm-penalty on portfolio weights. We follow the work of \cite{Ho2015} and consider a convex combinations of general $L_1$ and $L_2$-norm penalty functions, presented below.
\begin{equation}\label{eq:penalty_1}
P(\bz) = \alpha \bgamma_1 \lVert \bE(\btheta_1)\bz \rVert_1 + (1 - \alpha) \frac{\bgamma_2}{2} \lVert \bD(\btheta_2) \bz \rVert_2^2.
\end{equation}
Here, $\bE(\btheta_1)$ and $\bD(\btheta_2)$ are parameterized sparsifying and regularization transforms and therefore Equation $\eqref{eq:penalty_1}$ generalizes the standard elastic net formulation \citep{Zou2005}. The choice of the generalized elastic net penalty is deliberate and motivated by prior work that establishes the connection between norm-penalized portfolio weights with covariance shrinkage and robust optimization \citep{Demiguel2009}.

Determining the penalty parameters, $(\btheta_1,\btheta_2,\bgamma_1,\bgamma_2)$,  is generally difficult. Traditionally, the penalty structures, $\bE(\btheta_1)$ and $\bD(\btheta_2)$, are specified in advanced by the user and the amount of penalization,  $(\bgamma_1,\bgamma_2)$, is then estimated by cross-validation \citep{Demiguel2009, Elma2020, Hastie2001}. Recently, \cite{Butler2021IPOb, Butler2020IPO} present a parameter estimation approach that optimizes prediction model parameters to directly minimize downstream portfolio objectives. Indeed, more generally there is a growing body of research on the efficacy of integrated prediction and decision optimization modelling (see for example \cite{Amos2019, Donti2017, Bert2020, Elma2020, Elma2020b, Grigas2021, Mandi2020}). 

Motivated by the aforementioned work, in this paper we present an end-to-end framework for optimizing the penalty parameters of Equation $\eqref{eq:penalty_1}$ such that the resulting norm-penalized MVO program induces optimal portfolio decision-making. We structure the parameter estimation problem as  an end-to-end neural network with differentiable QP layers and present a novel technique for computing the derivative of the optimal primal solution with respect to the parameterized $L_1$-norm penalty by implicit differentiation of the dual program.  An advantage of our approach is that the penalty structure does need not be fully-specified in advance and instead penalty parameters are optimized in a data-driven manner. To our knowledge, our framework is the first to consider the optimization of parameterized norm-penalty structures in a portfolio optimization setting.  The remainder of the paper is outlined as follows. We begin with a review of the relevant literature in the field of robust and norm-penalized portfolio optimization. In Section \ref{sec:method} we present the optimization framework for learning  parameterized quadratic programs (QPs).  We demonstrate that the dual of the norm-penalized QP is a box-constrained QP and provide a first-order gradient descent approach for optimizing  the penalty parameters by backpropagating through differentiable QP layers \citep{Amos2017, Butler2021ADMM}. We conclude with a simulation study using both US stocks and global futures data and demonstrate the flexibility and effectiveness of the approach.  We find that the  parameterized norm-penalized MVO portfolios, in which the exact penalty structure is learned from the data, result in improved out-of-sample portfolio objectives.

\subsection{Relevant literature:}\label{sec:lit}
MVO portfolios are known to be sensitive to estimation error and thus frameworks for mitigating the adverse effects of estimation error are fundamental to the successful application of modern portfolio theory \citep{Chopra1993, Jobson1982, Michaud2008b}. For example,  \cite{Black1991} present a Bayesian framework for altering MVO inputs based on prior estimates of means and covariances. Alternatively, \cite{Michaud2008a} propose a resampled MVO process that averages the optimal weights derived from multiple resampled realizations of asset returns. However, the estimates assume a multivariate normal distribution, with moments equal to the in-sample estimates and the impact of sampling error is assumed to be in proportion to the variance of the sample mean. 

\cite{Ledoit2004,Ledoit2012}, demonstrate that when the covariance matrix is poorly conditioned the portfolio construction process can become unstable, resulting in concentrated `error-maximizing' portfolio weights. Their methods focus on providing improved covariance estimates based on James-Stein shrinkage \citep{Stein1956}. They propose covariance estimates that are a weighted average of the sample covariance estimate and a structured covariance matrix. The `optimal' amount of shrinkage is then determined by minimizing the mean squared error of the estimate to the realized covariance. 

Most closely related to our work, \cite{Demiguel2009} and \cite{Ho2015} mitigate `error-maximization' by augmenting MVO portfolios with $L_1$ and $L_2$-norm penalties on portfolio weights. They demonstrate that their norm-penalty formulations are equivalent to a robust MVO counterpart and propose tuning the penalty parameters by bootstrapping the historical distribution. More generally, robust optimization insulates against estimation error by minimizing objectives under worst-case realizations of the estimates \citep{Goldfarb2003}. The success of robust optimization typical depends on properly defining and calibrating the uncertainty structure and in general poorly calibrated uncertainty structures will result in overly conservative or intractable programs \citep{Bert2014}.  

While the seminal work in the field of robust optimization provides several theoretical methods for constructing good uncertainty sets  \citep{Ben2009,Ben2000,Bert2004,Goldfarb2003}, in quantitative finance, defining appropriate uncertainty structures remains a relatively ad-hoc process. For example, \cite{Goldfarb2003} propose calibrating the uncertainty structure for both the mean and covariance by considering the covariances estimates from a linear regression model and cast the robust counterpart as a convex second-order cone program.  \cite{Tu2004} consider a robust MVO counterpart with box uncertainty sets and calibrate the uncertainty bounds by bootstrapping the historical return distribution. \cite{Zhu2009} consider robust MVO under ellipsoidal uncertainty and demonstrate that when the uncertainty matrix is proportional to the asset covariance matrix then the robust MVO counterpart is equivalent to a nominal MVO with a larger risk aversion parameter. More recently, \cite{Yin2019} argue for the use of ellipsoidal uncertainty sets instead of the more restrictive box uncertainty sets  and propose calibrating the level of uncertainty in proportion to asset Sharpe ratios. 

In all cases discussed above we find that the regularization structures are determined in advance by the user and then the desired amount of regularization is estimated independently from the downstream portfolio optimization. In this paper, the portfolio regularization structure takes the form of parameterized norm-penalties on portfolio weights. In contrast to prior work, however, we present an end-to-end framework for optimizing parameterized norm-penalties in order to directly minimize the downstream portfolio optimization objective. Our parameterized norm-penalty optimization solution is motivated by recent work that advocates for the use of neural networks to optimize model hyper-parameters in settings where an exhaustive parameter search would be computationally impractical (see for example \cite{Feng2017,  Ich2021, Lorraine2018, Ped2016}). The proposed framework is flexible and obviates the need to pre-define a specific norm-penalty structure and instead enables users to formulate parameterized models that are fit to a  training dataset in a systematic and data-driven manner. Indeed,  previous work considers the integration of prediction modelling with downstream portfolio optimization using an end-to-end neural network approach \citep{Butler2021IPOb, Butler2020IPO,Chev2022,Uysal2021}. In all cases, however, the optimal portfolios are determined by differentiable optimization layers that directly compute the optimal primal solutions. However, forward solving and backpropagating through an $L_1$-norm penalized QP is challenging and is currently not supported by existing QP layer architectures. In this paper we present an alternative neural network architecture for  $L_1$-norm penalized QPs that instead solves the dual program - a box-constrained QP - and then recovers the primal solution from the optimal dual variables. To our knowledge, our end-to-end neural network is the first to consider implicit differentiation of the dual of an $L_1$-norm penalized QP and thus advances the understanding and implementation of differentiable optimization layers for solving non-smooth optimization programs.

\section{Methodology} \label{sec:method}
We follow closely to the methodology as described in \cite{Butler2020IPO}. We denote the matrix of (excess) return observations as ${\bY=[\by^{(1)},\by^{(2)},...,\by^{(m)} ] \in \mathbb{R}^{d_z \times m}}$  and let $\bz \in \mathbb{R}^{d_z}$ denote the vector of portfolio weights. The mean-variance optimal portfolio at time $i$ is given by:
\begin{equation}\label{eq:mvo_ineq}
\begin{split}
\minimize_z \quad &  -\bz^T \toi{ \by } + \frac{\delta}{2} \bz^T \toi{\bV} \bz \\
\text{subject to} \quad &   \bA \bz = \blb, \qquad \bG \bz \leq \bh\\
\end{split}
\end{equation}
with risk-aversion parameter $\delta \in \mathbb{R}_+$ and time-varying symmetric positive definite covariance matrix  $\toi{ \bV } \in \mathbb{R}^{d_z \times d_z}$. The linear equality and inequality constraints are described by matrices $\bA \in \mathbb{R}^{d_{\text{eq}} \times d_z }$, $\blb \in \mathbb{R}^{d_{\text{eq}}}$ and $\bG \in \mathbb{R}^{d_{\text{iq}} \times d_z}$, $\bh \in \mathbb{R}^{d_{\text{iq}}}$, respectively.

At decision time, $i-1$, however, the realization of asset returns, $\toi{ \by }$, is unobservable. Instead, associated feature variables, $\toi {\bx } \in \mathbb{R}^{d_x}$, are used to estimate $\toi{ \by }$. We consider linear regression models of the form:
$$\toi {\hat{\by} } =  \hat{\bbeta}^T \toi {\bx} $$
where $\hat{\bbeta} \in \mathbb{R}^{d_x \times d_y}$ denotes the matrix of regression coefficients. Let ${\ell(\hat{\by},\by) = \lVert  \hat{\by} - \by  \rVert_2^2}$ denote the least-squares loss function. Therefore, given a training data set $\mathcal{D} = \{(  \toi{ \bx },  \toi{ \by } )\}_{i=1}^m$, we estimate $\hat{\bbeta}$ by solving:
\begin{equation} \label{beta_hat}
\begin{split}
    \hat{\bbeta} & = \argmin_{  \bbeta }  \mathbb{E}_D[\ell (\bbeta^T \toi {\bx},  \toi{ \by } )],
\end{split}
\end{equation}
where  $\mathbb{E}_\mathcal{D}$ denotes the expectation with respect to the training distribution $\mathcal{D}$.  Following \cite{Butler2021IPOb}, the least-squares estimates of asset mean returns and covariances are given by:

\begin{equation} \label{eq:y_hat}
\begin{split}
\toi{\hat{  \by}} & =  \hat{\bbeta}^T  \toi{ \bx }\\
\toi{\hat{ \bV } } & = \hat{\bbeta}^T  \toi{\hat{\bW}} \hat{\bbeta} + \toi{\hat{\bF}},
\end{split}
\end{equation}
where $\toi{\hat{\bW }}$ denotes the time-varying feature covariance matrix and  $\toi{\hat{\bF}}$  denotes the diagonal matrix of residual variances. 

Using the estimates, $\toi{\hat{  \by}}$ and $\toi{\hat{ \bV } } $, we can solve the nominal MVO portfolio using information available at time $i-1$:
\begin{equation}\label{eq:mvo_ineq_est}
\begin{split}
\minimize_z \quad &   -\bz^T \toi{\hat{\by}} + \frac{\delta}{2} \bz^T \toi{\hat{ \bV } }  \bz \\
\text{subject to} \quad &   \bA \bz = \blb, \qquad \bG \bz \leq \bh.\\
\end{split}
\end{equation}
We denote the optimal MVO portfolio weights as $\toi{\bz^*}$. We reiterate that Program  $\eqref{eq:mvo_ineq}$ assumes that asset means and covariances are known with certainty, whereas Program $\eqref{eq:mvo_ineq_est}$ is formed at time $i-1$ using the estimates of mean returns and covariances. 

In reality all prediction models do make some error  (i.e. $ \toi{\hat{  \by}} \neq \toi{ \by} $, $\toi{\hat{\bV}} \neq \toi{\bV}$). In this paper we address the adverse effects of estimation error by augmenting the nominal MVO program with norm-penalties on portfolio weights:
\begin{equation}\label{eq:mvo_ineq_ppqp}
\begin{split}
\minimize_z \quad &  c_P(\bz) =   -\bz^T \toi{\hat{\by}} + \frac{\delta}{2} \bz^T \toi{\bV} \bz + \alpha \bgamma_1 \lVert \bE(\btheta_1)\bz \rVert_1 + (1 - \alpha) \frac{\bgamma_2}{2} \lVert \bD(\btheta_2)\bz \rVert_2^2\\
\text{subject to} \quad &   \bA \bz = \blb, \quad \bG \bz \leq \bh.
\end{split}
\end{equation}
Here $\lVert \cdot \rVert_p$ denotes the standard $L_p$-norm and $ 0 \leq \alpha \leq 1$. The matrices, $\bE(\btheta_1) \in \mathbb{R}^{d_{\text{E}} \times d_z}$  and ${\bD(\btheta_2) \in \mathbb{R}^{d_{\text{D}} \times d_z}}$ are parameterized sparsifying and regularization transforms, respectively.  Observe, that when $\alpha = 1$ and $\bE(\btheta_1) = \bI$ then the penalty is the standard $L_1$-norm penalty common to many applications in statistics and engineering \citep{Tibshirani1996, Kim2008}. Conversely, when $\alpha = 0$ and $\bD(\btheta_2) = \bI$ then the penalty is the standard $L_2$-norm penalty or Tikhonov regularization for approximations to ill-posed problems \citep{Tikhonov1963}. In its general form the norm-penalty resembles a generalization of the elastic net penalty  for joint regularization and feature selection in least-squares regression \citep{Zou2005}. As mentioned earlier, the choice of this specific norm-penalty is deliberate and motivated by prior work that establishes the connection between norm-penalized portfolio weights with covariance shrinkage and robust optimization.

The primary technical challenge is in optimizing the norm-penalty parameters such that the norm-penalized MVO portfolio induces optimal portfolio decisions. As described in Section \ref{sec:lit}, a traditional approach would pre-specify the norm-penalty structures, $\bE(\btheta_1)$ and $\bD(\btheta_2)$, and the optimal amount of penalization,  $( \bgamma_1, \bgamma_2)$, is then estimated independently from its use in the downstream portfolio optimization. In contrast, in the following subsection we present an end-to-end framework for optimizing parameterized norm-penalties with the aim of inducing optimal MVO portfolio decision-making.


\subsection{Optimization framework}\label{sec:method_sof}
Let $c \colon \mathbb{R}^{d_z \times m} \times  \mathbb{R}^{d_z \times m} \rightarrow \mathbb{R}$ denote the realized portfolio objective (cost) function, which takes as input the optimal decisions: $ \{ \toi{\bz^*} \}_{i=1}^m$ and realized returns: $\{ {\toi{\by}} \}_{i=1}^m$.  Specifically, in this paper we consider the minimum-variance, $c_{\text{MV}}$, and maximum Sharpe ratio, $c_{\text{SR}}$, costs, presented below:
\begin{equation} \label{eq:costs}
c_{\text{MV}}(\{ \toi{\bz^*} \}_{i=1}^m, \{ {\toi{\by}} \}_{i=1}^m )    =   \sigma^2  \quad \text{and} \quad  c_{\text{SR}}(\{ \toi{\bz^*} \}_{i=1}^m, \{ {\toi{\by}} \}_{i=1}^m )  = -\frac{ \mu-r_f }{ \sigma }
\end{equation}
where $r_f$ denotes the risk-free rate and:
$$ \mu = \frac{1}{m}\sum_{i=1}^m \toit{\bz^*} \toi{ \by } \quad \text{and} \quad \sigma^2 = \frac{1}{m} \sum_{i=1}^m(\toit{\bz^*} \toi{ \by } - \mu )^2.$$

We seek to optimize the penalty parameters, $\btheta_1$,$\btheta_2$, $\bgamma_1$ and $\bgamma_2$, in order to minimize the average realized cost, $c$, induced by the norm-penalized MVO portfolios $ \{ \toi{\bz^*} \}_{i=1}^m$.  We assume a discrete dataset $ \mathcal{D} = \{( \bx^{(i)}, \by^{(i)})\}_{i=1}^m$ and follow the work of \cite{Butler2021IPOb, Butler2020IPO} and cast the parameter estimation problem as a bi-level optimization program:
\begin{equation} \label{eq:stoch_discrete}
\begin{split}
\minimize_{ \btheta_1,\btheta_2, \bgamma_1, \bgamma_2 }  & \quad c(\{ \toi{\bz^*} \}_{i=1}^m, \{ {\toi{\by}} \}_{i=1}^m )  \\
\text{subject to } & \quad   \toi{\bz^*}  = \argmin_{\bz \in \mathbb{S}}  c_P( \bz ) \quad \forall \ i \in 1,...,m\\
& \quad \bgamma_1 \geq 0, \bgamma_2 \geq 0
\end{split}
\end{equation}
We denote the feasible region as $\mathbb{S} = \{ \bz \in \mathbb{R}^{d_z} | \bA \bz = \blb, \bG \bz \leq \bh \}$ and the objective function, $c_P$, is as defined in Program $\eqref{eq:mvo_ineq_ppqp}$. Observe that Program $\eqref{eq:stoch_discrete}$ results in a complicated dependency of the penalty parameters on the optimized value, $\bz^*$, connected through the $\argmin$ function.  In this paper we approximate a local solution to Program $\eqref{eq:stoch_discrete}$ by applying first-order gradient descent. Consider, for example, the gradient update $\partial  c/\partial  \btheta_1$; then by the chain-rule we have: ${ \partial  c/\partial  \btheta_1 = \partial  c/ \partial  \bz^* \cdot  \partial  \bz^*  /\partial \btheta_1}$. While computing the gradient, $ \partial  c/ \partial  \bz^* $, is relatively straightforward, computing the Jacobian, $\partial  \bz^*  /\partial \btheta_1$, requires differentiation through the $\argmin$ operator. To overcome this challenge, we make use of recent advances in neural network architecture for differentiable QPs and compute the action of the Jacobian by implicit differentiation of the optimality conditions of Program $\eqref{eq:mvo_ineq_ppqp}$  \citep{Amos2017, Agrawal2019, Butler2021ADMM}. This is described in more detail for general QPs in the following subsection.

\subsection{Differentiable QP layers} \label{sec:diff_qp}
We begin by first simplifying with respect to the $L_2$-norm penalty. Note that for ease of notation we temporarily drop the time index $i$.  We let:
\begin{equation}\label{eq:q_l2}
\hat{\bV}_{\bgamma_2} = \delta \hat{\bV} + (1 - \alpha) \bgamma_2 \bD(\btheta)^T\bD(\btheta)
\end{equation}
and therefore the norm-penalized objective of Program $\eqref{eq:mvo_ineq_ppqp}$ can be simplified as follows:
\begin{equation}\label{eq:ppqp_obj}
\begin{split}
c_P(\bz) = -\bz^T \hat{\by} + \frac{1}{2} \bz^T \hat{\bV}_{\bgamma_2}  \bz + \alpha \bgamma_1 \lVert \bE(\btheta)\bz \rVert_1.
\end{split}
\end{equation}
Observe that if $\alpha = 0$ then Program $\eqref{eq:mvo_ineq_ppqp}$ is a standard quadratic program and the relevant gradients with respect to the $L_2$-norm penalty are given as:
\begin{equation}\label{eq:grad_g2}
\begin{split}
\frac{\partial c   }{\partial \bgamma_2}  = (1 - \alpha) \Tr \Bigg( \frac{\partial c  }{\partial \hat{\bV}_{\bgamma_2}} \bD(\btheta)^T\bD(\btheta) \Bigg) \ \ \ \ & \frac{\partial c   }{\partial \bD} = \bgamma_2 (1 - \alpha) \Bigg(   \bD(\btheta) \frac{\partial c   }{\partial \hat{\bV}_{\bgamma_2}}  +  \bD(\btheta) \Big( \frac{\partial c   }{\partial \hat{\bV}_{\bgamma_2}}  \Big)^T \Bigg) 
\end{split}
\end{equation}
which ultimately requires computing the partial derivative $\partial c /  \partial \hat{\bV}_{\bgamma_2}$. In this paper, the gradients with respect to all QP program input variables are computed by restructuring program $\eqref{eq:stoch_discrete}$ as an end-to-end neural network with differentiable QP layers, described below.

The OptNet layer, presented by \cite{Amos2017}, is a differentiable optimization layer that uses an interior-point solver for forward solving batch quadratic programs. Backward differentiating is performed by implicit differentiation of a fixed-point provided by the KKT optimality conditions:

\begin{equation} \label{eq:dol_qp_kkt_g}
\begin{bmatrix}
-\hat{\by} +  \hat{\bV}_{\bgamma_2}\bz^*  + \bG^T \blambda^* + \bA^T \betta^* \\
\diag(\blambda^*) (\bG \bz^* - \bh)\\
\bA\bz^* - \blb
\end{bmatrix}
=
\begin{bmatrix}
\bzero \\
\bzero \\
\bzero
\end{bmatrix},
\end{equation}
where  $(\bz^*, \blambda^*, \betta^*)$ denotes the primal-dual solution of Program $\eqref{eq:mvo_ineq_ppqp}$. Taking the differential of  $\eqref{eq:dol_qp_kkt_g}$ and collecting like terms gives the following system of equations:
\begin{equation} \label{eq:dol_qp_diff}
\begin{split}
\begin{bmatrix}
  \hat{\bV}_{\bgamma_2}&  \bG^T &   \bA^T \\
\diag(\blambda^*)\bG & \diag (\bG \bz^* - \bh ) & 0\\
\bA & 0 & 0
\end{bmatrix}
\begin{bmatrix}
{\bdz }\\
{\bdlambda }\\
{\bdetta}
\end{bmatrix}
= -
\begin{bmatrix}
{{\bf d}\hat{\bV}_{\bgamma_2}} \bz^* - {{\bm d} \hat{\by} } + {\bdG}^T \blambda^* + {\bdA}^T \betta^* \\
\diag(\blambda ^*) {\bdG } \bz^* - \diag(\blambda ^*) {\bdh}\\
{\bdA} \bz^* - {\bdb}
\end{bmatrix}.
\end{split}
\end{equation}
As an example, the gradient, $\partial c / \partial \hat{\by}$ is then computed by the chain-rule: 
$$\frac{\partial c}{\partial \hat{\by}} = \frac{\partial c}{\partial \bz^*} \frac{\partial \bz^*}{\partial \hat{\by}}.$$
However, in the backpropagation algorithm it is inefficient to explicitly form the right-side Jacobian matrix and instead we compute  $\partial c /\partial \bz^*$ directly by solving the following system of equations: 

\begin{equation} \label{eq:dol_qp_sys}
\begin{split}
\begin{bmatrix}
\bar{ \bd }_{\bz} \\
\bar{ \bd }_{\blambda} \\
\bar{ \bd }_{\betta}
\end{bmatrix}
= -
\begin{bmatrix}
 \hat{\bV}_{\bgamma_2}&  \bG^T\diag(\blambda^*) &   \bA^T \\
\bG & \diag (\bG \bz^* - \bh ) & 0\\
\bA & 0 & 0
\end{bmatrix} ^{-1}
\begin{bmatrix}
\big( \frac{\partial c }{\partial \bz^*} \big)^T \\
0\\
0
\end{bmatrix}.
\end{split}
\end{equation}
Finally, the gradients of the cost with respect to all QP input variables are presented below and we refer to \cite{Amos2017} for more detail.

\begin{equation}\label{eq:dol_qp_grad}
\begin{aligned}
\frac{\partial c   }{\partial \hat{\bV}_{\bgamma_2}} & = \frac{1}{2} \Big(\bar{ \bd }_{\bz}  \bz^{*T} + \bz^* \bar{ \bd }_{\bz}^T \Big) & \qquad \frac{\partial c   }{\partial \hat{\by}} & = - \bar{ \bd }_{\bz} \\
\frac{\partial c   }{\partial \bA} & =  \bar{ \bd }_{\betta}  \bz^{*T} + \betta^* \bar{ \bd }_{\bz} ^T   & \qquad \frac{\partial c   }{\partial \blb} & = -\bar{ \bd }_{\betta}  \\
\frac{\partial c   }{\partial \bG} & = \diag(\blambda ^*)  \bar{ \bd }_{\blambda}  \bz^{*T} + \blambda^* \bar{ \bd }_{\bz} ^T    & \qquad  \frac{\partial c   }{\partial \bh} & = - \diag(\blambda ^*) \bar{ \bd }_{\blambda}
\end{aligned}
\end{equation}

More recently, \cite{Butler2021ADMM} provide an efficient differentiable optimization layer for forward solving and backward differentiating box-constrained QPs. Their framework applies the ADMM algorithm to solve the QP in the forward pass. The ADMM iterations are presented below:

\begin{subequations} \label{eq:admm_qp_iter_simp}
\begin{align}
\begin{bmatrix}
\tilde{\bz}^{k+1}\\
\betta^{k+1}
\end{bmatrix}
& = -
\begin{bmatrix}
 \hat{\bV}_{\bgamma_2}+ \rho \bI_{\bz} &   \bA^T \\
\bA & 0
\end{bmatrix}^{-1}
\begin{bmatrix}
 -\hat{\by} - \rho (\bz^{k} - \bmu^{k})\\
 -\blb
\end{bmatrix} \label{eq:admm_qp_x_iter_simp}\\
\bz^{k+1} &  = \Pi( \tilde{\bz}^{k+1}  + {\bmu}^k ) \label{eq:admm_qp_z_iter_simp}\\
{\bmu}^{k+1} & = {\bmu}^{k} +  \tilde{\bz}^{k+1} - \bz^{k+1} \label{eq:admm_qp_mu_iter_simp}
\end{align}
\end{subequations}
where $\rho > 0 $ is a user-defined penalty parameter and $\Pi$ denotes the euclidean projection onto a set of box constraints $\{ \bz \in \mathbb{R}^{d_z} \mid \bl \leq  \bz \leq \bu \}$. The dual variables associated with the box constraints are defined as $\tilde{\blambda}^* = (\blambda^*_-,\blambda^*_+)$ with:
\begin{equation} \label{eq:lambda_admm}
\blambda^*_- = -\min(\rho \bmu^*,0) \quad \text{and} \quad  \blambda^*_+ = \max(\rho \bmu^*,0).
\end{equation}
The gradients of the optimal solution, $\bz^*$, with respect to the QP input variables are computed by implicit differentiation of a custom fixed-point mapping. Let $\hat{ \bd }_{\bz}$ and $\hat{ \bd }_{\betta}$ be defined as:
\begin{equation}\label{eq:grads_admm}
\begin{split}
\begin{bmatrix}
\hat{ \bd }_{\bz}  \\
\hat{ \bd }_{\betta}
\end{bmatrix}
& =
 \Bigg[
\begin{bmatrix}
D\Pi(\bp) & 0\\
0 & \bI_{\betta}
\end{bmatrix}
\begin{bmatrix}
 \hat{\bV}_{\bgamma_2}+ \rho \bI_{\bz} &   \bA^T \\
\bA & 0
\end{bmatrix}
+
\begin{bmatrix}
- \rho (2D\Pi(\bp) - \bI_{\bz}) & 0\\
0 & 0
\end{bmatrix}
\Bigg]^{-1}
\begin{bmatrix}
D\Pi(\bp) & 0\\
0 & \bI_{\betta}
\end{bmatrix}
\begin{bmatrix}
\big( - \frac{\partial c }{\partial \bz^*} \big)^T \\
0
\end{bmatrix},
\end{split}
\end{equation}
where $D\Pi$ denotes the derivative of the projection operator and $\bp = \bz^* + \bmu^*$. Then the gradients of the cost function, $c$, with respect to problem variables $\hat{\bV}_{\bgamma_2}$, $\hat{\by}$, $\bA$ and $\blb$ are given by: 
\begin{equation}\label{eq:admm_partials}
\begin{aligned}
\frac{\partial c   }{\partial \hat{\bV}_{\bgamma_2}} & = \frac{1}{2} \Big(\hat{ \bd }_{\bz}   \bz^{*T} + \bz^* \hat{ \bd }_{\bz}^T \Big) & \qquad \frac{\partial c   }{\partial \hat{\by}} & = -\hat{ \bd }_{\bz}  \\
\frac{\partial c   }{\partial \bA} & =  \hat{ \bd }_{\betta}  \bz^{*T} + \betta^* \hat{ \bd }_{\bz} ^T   & \qquad \frac{\partial c   }{\partial \blb} & = -\hat{ \bd }_{\betta}
\end{aligned}
\end{equation}
Similarly,  we define $\tilde{\bmu}^*$ and $\hat{ \bd }_{\blambda}$ as:
\begin{equation}\label{eq:mu_tilde}
\tilde{\bmu}^*_j =  \begin{cases}
                \bmu^*_j & \text{if } \bmu^*_j \neq 0\\
                1 & \text{otherwise} \\
                \end{cases}
 \qquad \text{and} \qquad \hat{ \bd }_{\blambda} = \diag(\rho \tilde{\bmu}^*)^{-1} \Big( - \Big(\frac{\partial c }{\partial \bz^*} \Big)^T -  \hat{\bV}_{\bgamma_2}\hat{ \bd }_{\bz} - \bA^T \hat{ \bd }_{\betta} \Big).
\end{equation}
The gradients of the cost function, $c$, with respect to the box constraint variables are given as:
\begin{equation}\label{eq:admm_partials_box}
\begin{aligned}
\frac{\partial c   }{\partial \bl} & = \diag(\blambda^*_-)\hat{ \bd }_{\blambda} & \qquad  \frac{\partial c   }{\partial \bu} & = -\diag(\blambda^*_+)\hat{ \bd }_{\blambda}.
\end{aligned}
\end{equation}

\subsection{Implicit dual differentiation} \label{sec:dual}
In the general case, Program $\eqref{eq:mvo_ineq_ppqp}$ is a $L_1$-norm penalized QP. Note that the $L_1$-norm is not differentiable everywhere and is certainly not twice differentiable; thus complicating the implicit $\argmin$ differentiation. In this section we demonstrate that the dual of Program $\eqref{eq:mvo_ineq_ppqp}$ is a convex box-constrained QP. It is therefore possible to forward solve and backward differentiate through the dual program using the methodology described in Section \ref{sec:diff_qp}. The optimal portfolio, $\bz^*$, is then recovered from the optimal dual solution. 

 Let $\bw = \bE(\btheta)\bz$, then the following program is equivalent to Program $\eqref{eq:mvo_ineq_ppqp}$:
\begin{equation}\label{eq:ppqp_en_dual_start}
\begin{split}
\minimize_{\bz} \quad &   -\bz^T \hat{\by} + \frac{1}{2} \bz^T \hat{\bV}_{\bgamma_2}  \bz + \alpha \bgamma_1 \lVert \bw \rVert_1 \\
\text{subject to} \quad &   \bA  \bz = \blb, \quad \bG  \bz \leq \bh  \\
\quad & \bw = \bE(\btheta)\bz.
\end{split}
\end{equation}
The Lagrangian of Program $\eqref{eq:ppqp_en_dual_start}$  is presented below:
\begin{equation}\label{eq:ppqp_en_dual_l}
\begin{split}
\mathcal{L}(\bz, \bw, \bv, \betta, \blambda) = & -\bz^T \hat{\by} + \frac{1}{2} \bz^T \hat{\bV}_{\bgamma_2}  \bz + \alpha \bgamma_1 \lVert \bw \rVert_1 + \bv^T(\bE(\btheta)\bz - \bw )\\
 & + \betta^T(\bA  \bz - \blb) + \blambda^T(\bG  \bz - \bh )
 \end{split}
\end{equation}
with dual variables $\bv \in \mathbb{R}^{d_E}$, $\betta \in \mathbb{R}^{d_{eq}}$ and $\blambda \in \mathbb{R}_+^{d_{iq}}$. Following \cite{Kim2008}, the Lagrange dual problem associated with Program $\eqref{eq:ppqp_en_dual_start}$ is given by:

\begin{equation}\label{eq:ppqp_en_dual_problem}
\begin{split}
\maximize_{\bv, \betta, \blambda} \quad & \inf_{\bz,\bw} \mathcal{L}(\bz, \bw, \bv, \betta, \blambda)\\
\text{subject to} \quad &   \blambda \geq 0.
 \end{split}
\end{equation}
Observe that Program $\eqref{eq:ppqp_en_dual_start}$ is convex and therefore strong duality holds \citep{Boyd2004}.
The first-order optimality conditions with respect to the primal variable, $\bz$, results in the following linear system of equations:
\begin{equation}\label{eq:ppqp_en_dual_z}
\begin{split}
\hat{\bV}_{\bgamma_2} \bz = \hat{\by} - \bE(\btheta)^T\bv - \bA^T\betta - \bG^T\blambda.
 \end{split}
\end{equation}
Furthermore, dual feasibility is given by boundedness with respect to the dual variable, $\bv$, specifically:
\begin{equation}\label{eq:ppqp_en_dual_w}
\begin{split}
\inf_{\bw} \alpha \bgamma_1 \lVert \bw \rVert_1 - \bv^T\bw & = -\sup_{\bw} \bv^T\bw - \alpha \bgamma_1 \lVert \bw \rVert_1\\
& = \begin{cases} 0 & \text{if } \lVert \bv \rVert_\infty \leq \alpha \bgamma_1 \\ -\infty & \text{otherwise}. \end{cases}
 \end{split}
\end{equation}
Combining the results from Equations $\eqref{eq:ppqp_en_dual_z}$ and $\eqref{eq:ppqp_en_dual_w}$ results in the following box-constrained QP:
\begin{equation}\label{eq:ppqp_en_dual}
\begin{split}
\minimize_{\bv, \betta, \blambda} \quad &   \frac{1}{2} \br^T \hat{\bV}^{-1}_{\bgamma_2} \br + \betta^T \blb + \blambda^T \bh\\
\text{subject to} \quad & -\alpha \bgamma_1 \bone \leq  \bv  \leq \alpha \bgamma_1 \bone, \quad \blambda \geq 0.
\end{split}
\end{equation}
where $\br = \hat{\by} - \bE(\btheta)^T\bv - \bA^T\betta - \bG^T\blambda$ and $\bone \in \mathbb{R}^{d_E}$ is a vector of ones. Note that in the absence of constraints, Program $\eqref{eq:ppqp_en_dual}$ reduces to a convex box-constrained QP in the dual variable $\bv$:
\begin{equation}\label{eq:ppqp_en_dual_no_con}
\begin{split}
\minimize_{\bv} \quad &   \frac{1}{2} \bv^T \bE(\btheta) \hat{\bV}^{-1}_{\bgamma_2} \bE(\btheta)^T \bv - \bv^T \bE(\btheta) \bV^{-1}_{\bgamma_2} \hat{\by}\\
\text{subject to} \quad & -\alpha \bgamma_1 \bone \leq  \bv  \leq \alpha \bgamma_1 \bone.
\end{split}
\end{equation}
In the general case we denote the primal-dual variables as $\bnu = [\bv, \betta, \blambda]$ and define the vector $\bq$ and positive semidefinite matrix $\bQ$ as:

\begin{equation}\label{eq:big_q}
\bq =  \begin{bmatrix}
-\bE(\btheta) \hat{\bV}^{-1}_{\bgamma_2} \hat{\by}\\
-\bA \hat{\bV}^{-1}_{\bgamma_2} \hat{\by} + \blb \\
- \bG \hat{\bV}^{-1}_{\bgamma_2}\hat{\by} + \bh
\end{bmatrix}
\quad 
\bQ = \begin{bmatrix}
\bE(\btheta) \hat{\bV}^{-1}_{\bgamma_2} \bE(\btheta)^T & \bE(\btheta) \hat{\bV}^{-1}_{\bgamma_2} \bA^T & \bE(\btheta) \hat{\bV}^{-1}_{\bgamma_2} \bG^T\\
\bA \hat{\bV}^{-1}_{\bgamma_2} \bE(\btheta)^T		     &	\bA \hat{\bV}^{-1}_{\bgamma_2} \bA^T              & \bA \hat{\bV}^{-1}_{\bgamma_2} \bG^T\\
 \bG \hat{\bV}^{-1}_{\bgamma_2} \bE(\btheta)^T             &  \bG \hat{\bV}^{-1}_{\bgamma_2} \bA^T             &   \bG \hat{\bV}^{-1}_{\bgamma_2} \bG^T 
\end{bmatrix}.
\end{equation}
Then we have the following equivalent convex box-constrained QP:
\begin{equation}\label{eq:ppqp_en_dual_q}
\begin{split}
\minimize_{\bnu} \quad &   \frac{1}{2} \bnu^T \bQ  \bnu +\bnu^T \bq\\
\text{subject to} \quad & \bl \leq \bnu \leq \bu,
\end{split}
\end{equation}
with box constraints given as:
\begin{equation}\label{eq:box}
\bl = \begin{bmatrix}
-\alpha \bgamma_1 \bone\\
\bm{-\infty}\\
\bzero
\end{bmatrix} \qquad
\bu = \begin{bmatrix}
\alpha \bgamma_1 \bone \\
\bm{\infty}\\
\bm{\infty}
\end{bmatrix}
\end{equation}
We apply the methodology described in Section \ref{sec:diff_qp} for box-constrained QPs to forward solve and backward differentiate Program $\eqref{eq:ppqp_en_dual_q}$ with respect to the dual program variables: $\bQ$, $\bp$, $\bl$, and $\bu$. The optimal primal solution, $\bz^*$, is then determined by solving Equation $\eqref{eq:ppqp_en_dual_z}$ with respect to the optimal dual variables. Lastly, we note that implicit differentiation of the dual Program will often result in a near singular linear system of  Equations $\eqref{eq:grads_admm}$ and in practice we solve the $L_2$-regularized system by adding the regularization term $\epsilon \bI$ \citep{Tikhonov1963}. We find that regularization values of $\epsilon$ in the range of $ [ 10^{-8}, 10^{-4} ]$ work well for most applications.

\subsection{End-to-end neural network}\label{sec:network}
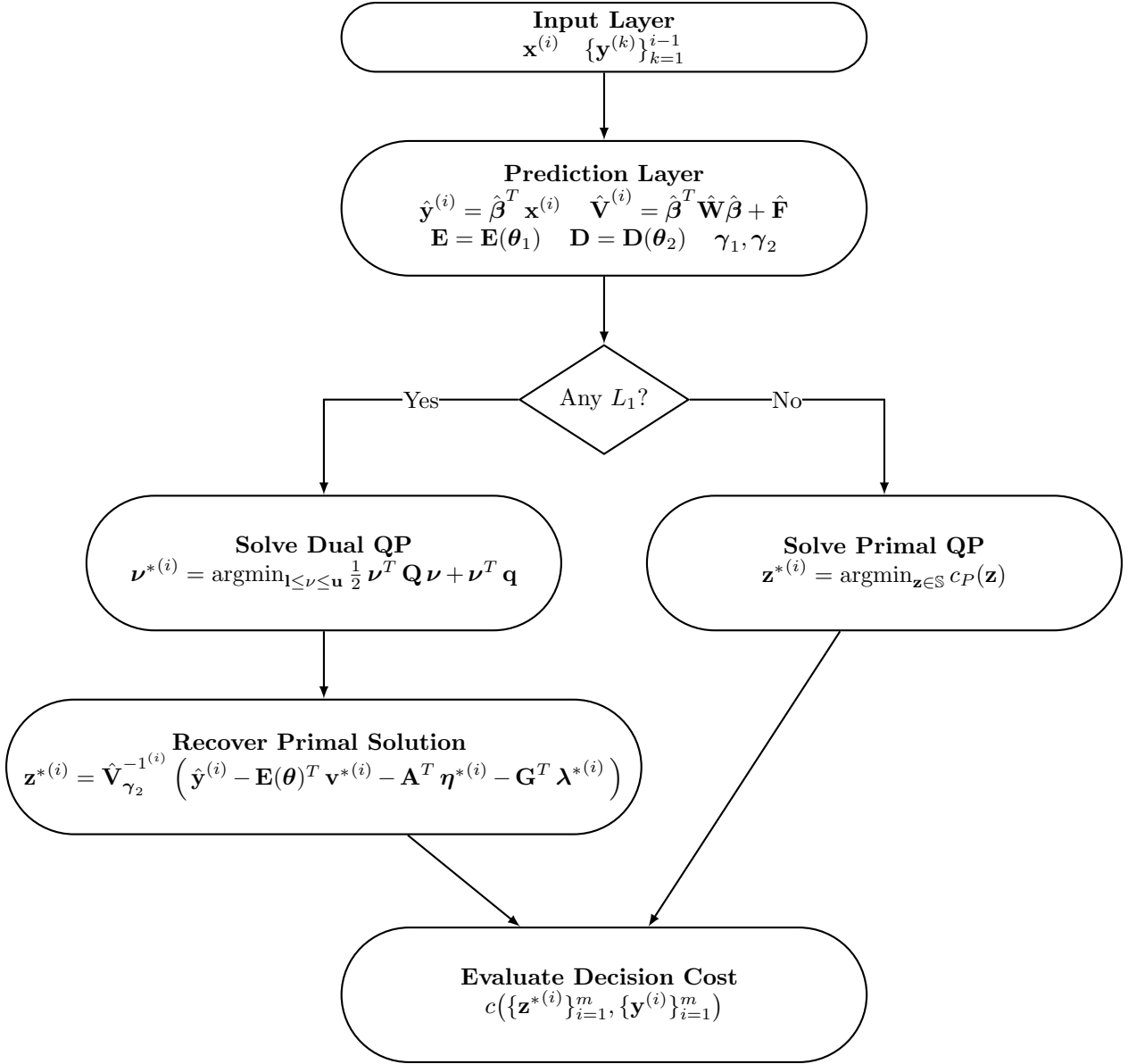
\begin{figure}[]
 \begin{tikzpicture}[font=\small,thick]
 
\node[draw,
    rounded rectangle,
    align = center,
    minimum width=8cm,
    minimum height=1cm] (block1) {{\textbf{Input Layer}} \\  $\toi{\bx} \quad \{ \toi[k]{\by} \}_{k = 1}^{i-1} $} ;
 
\node[draw,
rounded rectangle,
    align = center,
    minimum width=8cm,
    minimum height=2cm,
    below=of block1
] (block2){{\textbf{Prediction Layer}} \\  $\toi{\hat{  \by}}  =  \hat{\bbeta}^T  \toi{ \bx } \quad \toi{\hat{ \bV } }  = \hat{\bbeta}^T  \hat{\bW} \hat{\bbeta} + \hat{\bF} $ \\ 
$\bE = \bE(\btheta_1) \quad  \bD = \bD(\btheta_2) \quad \bgamma_1, \bgamma_2$
};

\node[draw,
    diamond,
    below=of block2,
    minimum width=2.5cm,
    inner sep=0] (block4) { Any $L_1$?};
 
\node[draw,
    rounded rectangle,
    align = center,
    below left=of block4,
     minimum width=7cm,
    minimum height=2cm,
    inner sep=0] (block5) { {\textbf{Solve Dual QP}} \\ $\toi{\bnu^*} = \argmin_{\bl \leq  \nu \leq \bu }  \frac{1}{2} \bnu^T \bQ  \bnu +\bnu^T \bq $};
 
\node[draw,
    rounded rectangle,
    align = center,
    below right=of block4,
     minimum width=7cm,
    minimum height=2cm,
    inner sep=0] (block6) { {\textbf{Solve Primal QP}} \\ $\toi{\bz^*} = \argmin_{\bz \in \mathbb{S}}  c_P(\bz) $};
 
\node[draw,
rounded rectangle,
    align = center,
    below=of block5,
    minimum width=7cm,
    minimum height=2cm] (block7) { {\textbf{Recover Primal Solution }} \\ $\toi{\bz^*} = \toinv{\hat{\bV}_{\bgamma_2}} \Big( \toi{\hat{\by}} - \bE(\btheta)^T\toi{\bv^*} - \bA^T\toi{\betta^*} - \bG^T\toi{\blambda^*} \Big)$};

 
 
 
\node[draw,
    rounded rectangle,
    align = center,
    below=7cm of block4,
    minimum width=8cm,
    minimum height=2cm,] (block11) { {\textbf{Evaluate Decision Cost }} \\ $c\big (\{ \toi{\bz^*} \}_{i=1}^m, \{ {\toi{\by}} \}_{i=1}^m \big)$ };
 
\node[coordinate,below=4.35cm of block4] (block12) {};

\draw[-latex] (block1) edge (block2)
    (block2) edge (block4)
    (block5) edge (block7)
    (block6) edge (block11)
    (block7) edge (block11)
    ;
 
\draw[-latex] (block4) -| (block5)
    node[pos=0.25,fill=white,inner sep=0]{Yes};
 
\draw[-latex] (block4) -| (block6)
    node[pos=0.25,fill=white,inner sep=0]{No};
 
 
 
 
\end{tikzpicture}
\caption{Bi-level optimization framework represented as an end-to-end neural network with prediction layer, differentiable quadratic programming layer and decision cost loss function. }
\label{fig:network}
\end{figure}

We follow the work of \cite{Donti2017, Amos2019, Mandi2019} and others, and solve Program $\eqref{eq:stoch_discrete}$ by  structuring the bi-level optimization program as an end-to-end trainable neural network with differentiable quadratic programming layers, depicted in Figure \ref{fig:network}. In the forward pass, the input layer takes the feature variables $\toi{\bx}$, and often also the historical realized returns $\{ \toi[k]{\by} \}_{k = 1}^{i-1}$, and passes them to a prediction model layer to produce the estimates: $\toi{\hat{\by} } $ and $\toi{\hat{ \bV } }$, as well as the norm-penalty structures, $\bE$ and $\bD$, and their respective magnitudes $\bgamma_1$ and $\bgamma_2$. The QP input variables are passed to a differentiable optimization layer which solves the norm-penalized MVO program and returns the optimal portfolio decisions: $\{ \toi{\bz^*} \}_{i=1}^m$. 
 
Note that if the penalized MVO program does not contain an $L_1$-norm penalty, then Program $\eqref{eq:mvo_ineq_ppqp}$ is a standard QP and can be solved by applying the differentiable QP layer(s) described in Section \ref{sec:diff_qp}. Otherwise, we solve Program $\eqref{eq:mvo_ineq_ppqp}$  by first solving the dual problem, using the methodology described in Section  \ref{sec:dual}. The primal solution is then recovered from the optimal dual variables. The quality of the optimal decisions are evaluated by the cost function, $c$, over the training dataset:  $\{\toi{\bx}, \toi{\by} \}_{i = 1}^m$. At each iteration, the penalty parameters are updated by applying the backpropagation algorithm \citep{Hinton1986}.

\subsection{Example}
The neural network framework depicted in Figure \ref{fig:network} is made available as a \textit{torch} module in the open-source R package \textit{lqp} (learning quadratic programs), available here: \textit{https://github.com/butl3ra/lqp}. We present an example of learning a parameterized norm-penalty of the form: 
$$
{P(\bz) =  \bgamma_1 \lVert \bE(\btheta_1)\bz \rVert_1 +  \frac{\bgamma_2}{2} \lVert \bD(\btheta_2) \bz \rVert_2^2},
$$
where $\bE(\btheta_1) = \Diag(\btheta_1)$, $\btheta_1 \in \mathbb{R}_+^{d_z}$, $\bD( \btheta_2) =  \hat{\bW}^{\frac{1}{2}}  \btheta_2 $, and $\hat{\bW}$ denotes the in-sample covariance matrix of feature variables. We consider a long-only, fully invested minimum-variance portfolio optimization on a universe of $50$ US stocks, with feature variables given by the Famma-French Five factor model. We refer the reader to the discussion in Section \ref{sec:results} for more comprehensive implementation details.

\begin{figure}[]
  \centering
  \begin{subfigure}[b]{0.45\linewidth}
    \includegraphics[width=\linewidth , trim={0mm 0cm 0cm 0cm},clip]{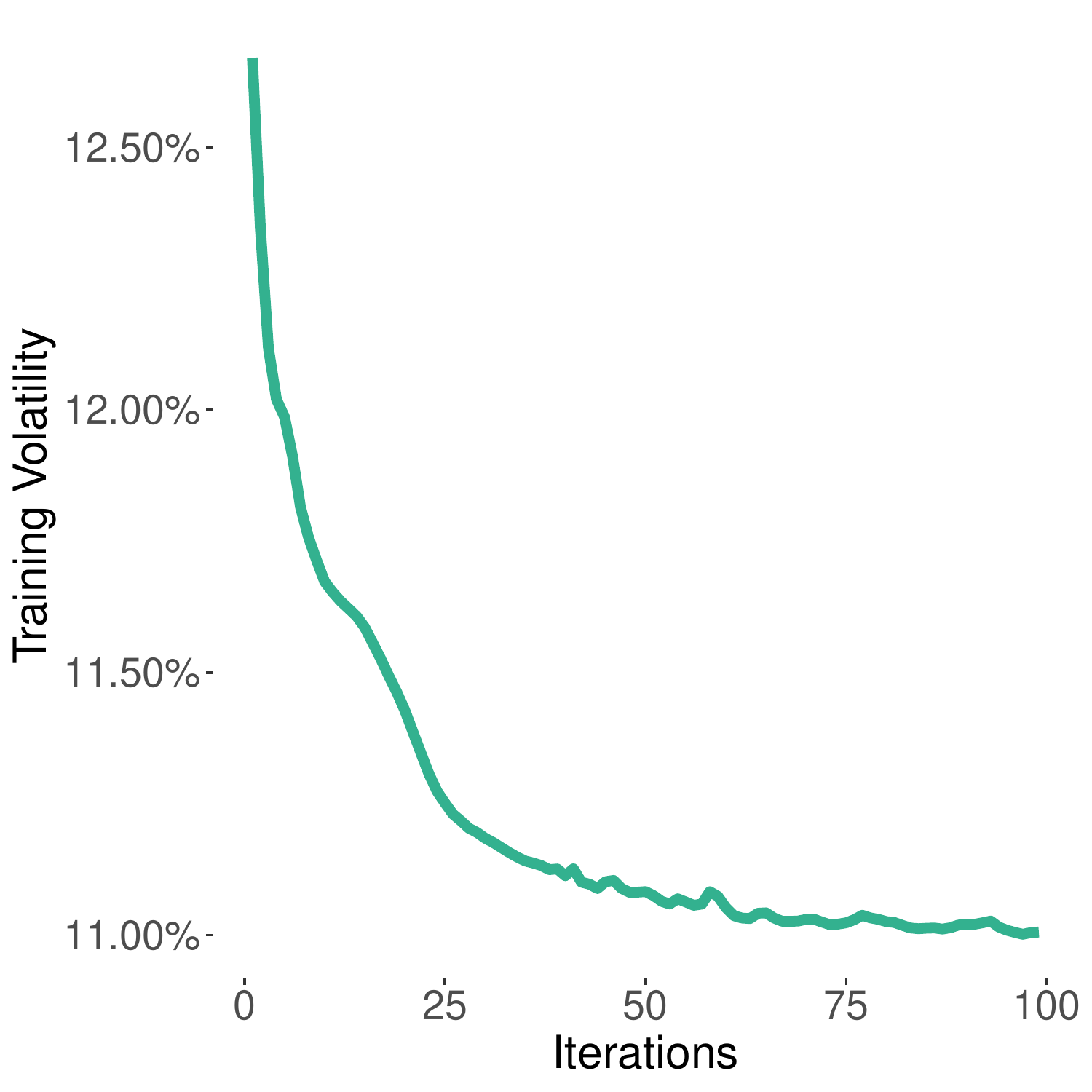}
    \caption{Training volatility.}
      \label{fig:demo_a}
  \end{subfigure}
  \begin{subfigure}[b]{0.45\linewidth}
    \includegraphics[width=\linewidth , trim={0mm 0cm 0cm 0cm},clip]{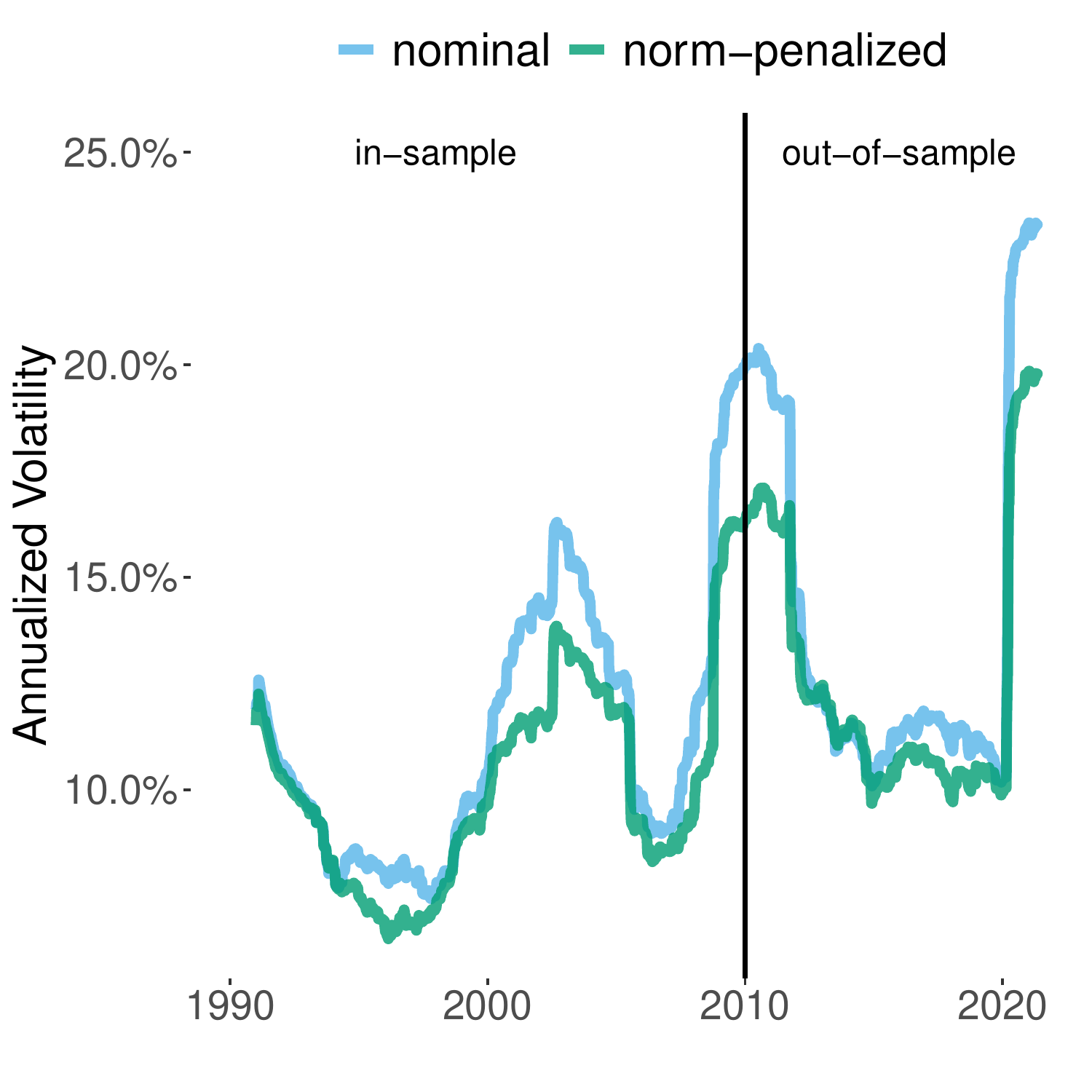}
    \caption{Rolling $3$-year volatility.}
      \label{fig:demo_b}
  \end{subfigure}
    \caption{In-sample training volatility and rolling $3$-year annualized volatility of the nominal and norm-penalized minimum-variance portfolios.}
  \label{fig:demo}
\end{figure}

Figure \ref{fig:demo_a} reports the in-sample training volatility of the parameterized norm-penalized minimum-variance portfolio as a function of the number of training iterations. Observe that the in-sample training volatility decreases from an annualized volatility of over $12.5\%$ to an annualized volatility of approximately $11.0\%$ after $100$ iterations of gradient descent. Figure \ref{fig:demo_b} charts the rolling $3$-year realized volatility of the nominal and norm-penalized portfolios. The in-sample training period is from January 1990 to December 2009 and the out-of-sample period is from January 2010 to May 2021. Observe that the norm-penalized portfolio provides a consistent reduction in realized volatility; in particular during periods of market stress, such as the Dot-com Bubble (2000-2003), the Global Financial Crisis (2008-2009) and the COVID-19 Pandemic (2020-2021).

\section{Experiments}\label{sec:results}
We present two experiments: a minimum-variance portfolio optimization on US stocks data and a maximum Sharpe ratio optimization on global commodity futures data.  The US stocks data  consists of $255$ liquid US stocks traded on major U.S. exchanges (NYSE, NASDAQ, AMEX, ARCA). Data is provided by Quandl and is summarized in Table \ref{table:stocks} in Appendix \ref{sec:app_data}. Weekly price data is given from January $1990$ through May $2021$. The futures markets are summarized in Table \ref{table:futures}. Futures data is given from March $1986$ through May $2021$, and is provided by Commodity Systems Inc. The returns of futures contracts are computed directly from the price data and are in excess of the risk-free rate.

We consider several norm-penalty models, summarized below.
\begin{enumerate}
\item \textbf{Nominal}: the nominal MVO program: $P(\bz) = 0.$
\item \textbf{L2}: user-defined $L_2$-norm penalty: $P(\bz) = \frac{\bgamma_2}{2} \lVert \bz \rVert_2^2.$
\item \textbf{L2-COV}: user-defined following \cite{Goldfarb2003} and letting $\bD = \hat{\bW}^{\frac{1}{2}} \hat{\bbeta}$, where $\hat{\bW}$ denotes the covariance matrix of feature variables: $P(\bz) = \frac{\bgamma_2}{2} \lVert \bD \bz \rVert_2^2.$
\item \textbf{L1}: user-defined $L_1$-norm penalty: $P(\bz) = \bgamma_1 \lVert \bz \rVert_1.$
\item \textbf{EN}: user-defined elastic-net  penalty with mixing coefficient $\alpha = 0.5$ and ${P(\bz) = \alpha \bgamma_1 \lVert \bz \rVert_1 + (1 - \alpha) \frac{\bgamma_2}{2} \lVert \bz \rVert_2^2 }.$
\item \textbf{L2-P}: parameterized $L_2$-norm penalty with $\bD(\btheta_2) = \diag(\btheta_2)$, $\btheta_2 \in \mathbb{R}_+^{d_z}$ and
${P(\bz) = \frac{\bgamma_2}{2} \lVert \bD(\btheta_2) \bz \rVert_2^2}.$
\item \textbf{L2-COV-P}: parameterized $L_2$-norm penalty with $\bD( \btheta_2) =   \hat{\bW}^{\frac{1}{2}}  \btheta_2 $ and ${P(\bz) = \frac{\bgamma_2}{2} \lVert \bD( \btheta_2) \bz \rVert_2^2}.$
\item \textbf{L1-P}: parameterized $L_1$-norm penalty with $\bE(\btheta_1) = \diag(\btheta_1)$, $\btheta_1 \in \mathbb{R}_+^{d_z}$ and
${P(\bz) = \bgamma_1 \lVert \bE(\btheta_1)\bz \rVert_1}.$
\item \textbf{EN-P}: parameterized elastic-net penalty with mixing coefficient $\alpha = 0.5$, ${\bD(\btheta_2) = \diag(\btheta_2)}$, $\bE(\btheta_1) = \diag(\btheta_1)$ and  ${P(\bz) = \alpha \bgamma_1 \lVert \bE(\btheta_1)\bz \rVert_1 + (1 - \alpha) \frac{\bgamma_2}{2} \lVert \bD(\btheta_2) \bz \rVert_2^2}.$
\end{enumerate}
We refer the reader to Appendix \ref{sec:app_imp_details} for full implementation details. 

Note that the norm-penalty models, described above, can be divided into two groups:
\begin{enumerate}
\item \textbf{User-defined:} the norm-penalty is pre-defined by the user and the amount of regularization is optimized by solving Program $\eqref{eq:stoch_discrete}$.
\item \textbf{Parameterized:} the  norm-penalty is parameterized and therefore the norm-penalty parameters and amount of regularization are jointly optimized  by solving Program $\eqref{eq:stoch_discrete}$.
\end{enumerate}

Both experiments consist of $30$ randomized trials where, at the beginning of each trial, a basket of $n$ assets is randomly drawn from the universe of available assets and held constant throughout the duration of the trial.  Portfolios are formed at the close of each week and rebalanced on a weekly basis. All experiments use data from January $1990$ through December $2009$ for training. Performance is evaluated over the out-of-sample period from January $2010$ through May $2021$.  Our analysis focuses on comparing the out-of-sample portfolio objectives, defined in Equation $\eqref{eq:costs}$, of the norm-penalized portfolios relative to the nominal MVO program. 

In all experiments we estimate the penalty parameters using the end-to-end neural network framework, described in Section \ref{sec:network}. Both the primal and dual norm-penalized QP layers are implemented by applying the ADMM layer, described in Sections \ref{sec:diff_qp} and \ref{sec:dual}. All portfolio optimizations under consideration are either unconstrained or contain linear equality constraints and box inequality constraints. Indeed, under these conditions, the ADMM layer implementation is preferred for its reported computational efficiency and ability to support larger scale portfolio optimization problems with $100 - 1000$ decision variables and thousands of training examples \citep{Butler2021ADMM}.

\subsection{Experiment 1: US stocks data}\label{sec:exp_1}
We consider the long-only, fully invested minimum-variance optimization $(\mathbb{S} = \{ \bz \mid \bone^T\bz = 1, \bz \geq 0 \} )$  on US stocks data, described in Table \ref{table:stocks}. The norm-penalty parameters are determined by solving Program $\eqref{eq:stoch_discrete}$  to minimize the  variance decision cost, $c_{\text{MV}}$. Asset mean returns are estimated by least-squares according to Equation $\eqref{eq:y_hat}$. The feature variables, $\toi{\bx} \in \mathbb{R}^5$, are given by the Famma-French Five factor model with data provided by the Kenneth R. French data library. We estimate the time-varying covariance of feature variables, $\toi{\hat{\bW}}$,  by applying a simple moving average covariance evaluated over the trailing $52$-weeks. The asset covariance matrix, $\toi{\hat{ \bV } } $, is given by  Equation $\eqref{eq:y_hat}$.  In each trial experiment we randomly draw $n=50$ stocks from the asset universe. 
\begin{figure}[]
  \centering
    \includegraphics[width=0.75\linewidth , trim={0mm 0cm 0cm 0cm},clip]{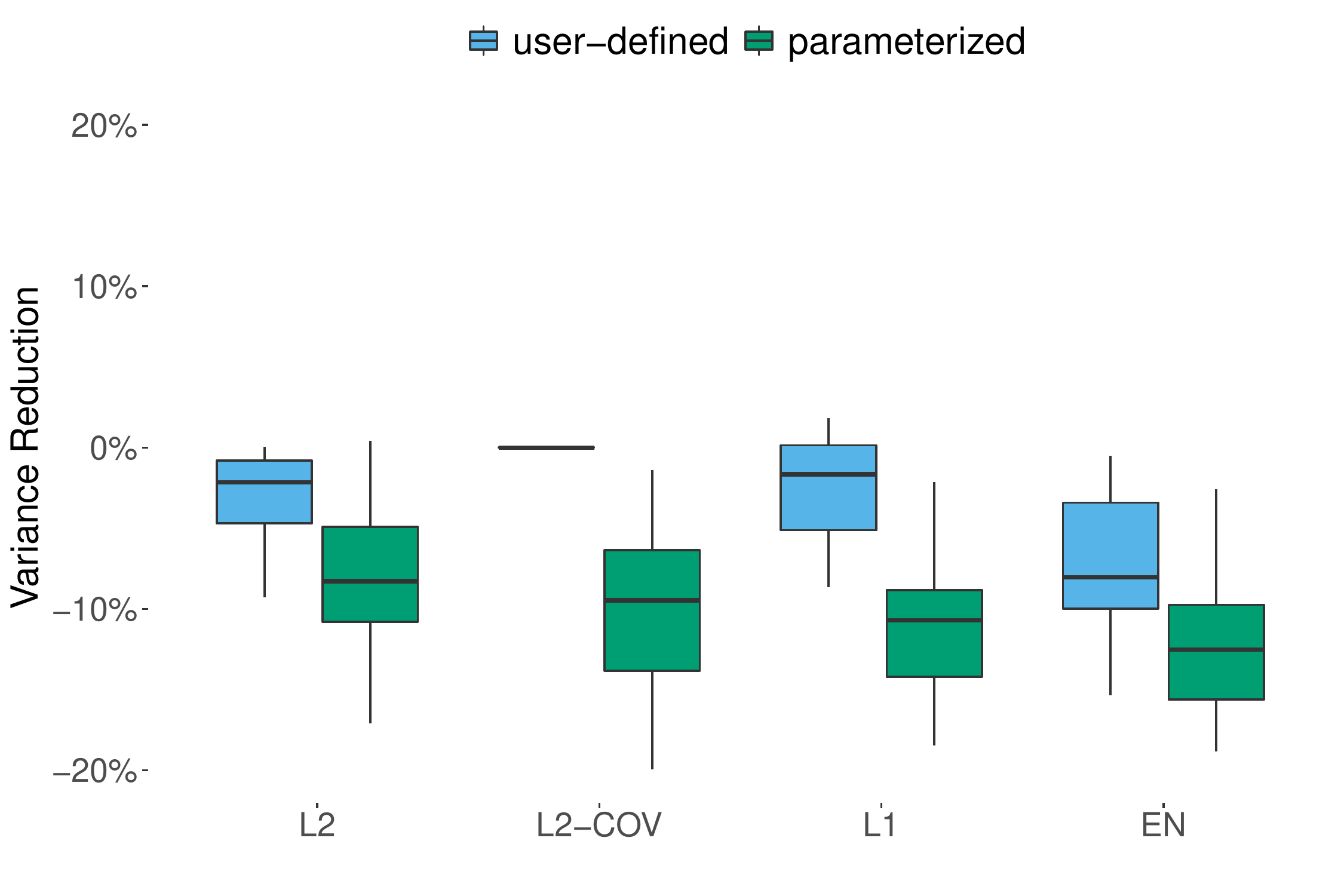}
   \caption{Out-of-sample percent variance cost reduction of norm-penalized minimum-variance portfolios relative to the nominal minimum-variance portfolio evaluated over $30$ randomized trials}
  \label{fig:mv_box}
\end{figure}

Figure \ref{fig:mv_box} reports the out-of-sample percent reduction in variance of the norm-penalized minimum-variance portfolios relative to the nominal minimum-variance portfolio. With the exception of the L2-COV model, all norm-penalized portfolios reduce the out-of-sample portfolio variance. We observe that the user-defined models produce a marginal $0\% - 7\%$ average reduction in out-of-sample variance. In contrast, the parameterized models reduce the out-of-sample variance by $9\%- 14\%$ on average. For both user-defined and parameterized models, the EN penalty, which contains both $L1$ and $L2$-norm penalties, produces the smallest average out-of-sample variance.

\begin{figure}[]
  \centering
  \begin{subfigure}[b]{0.45\linewidth}
    \includegraphics[width=\linewidth , trim={0mm 0cm 0cm 0cm},clip]{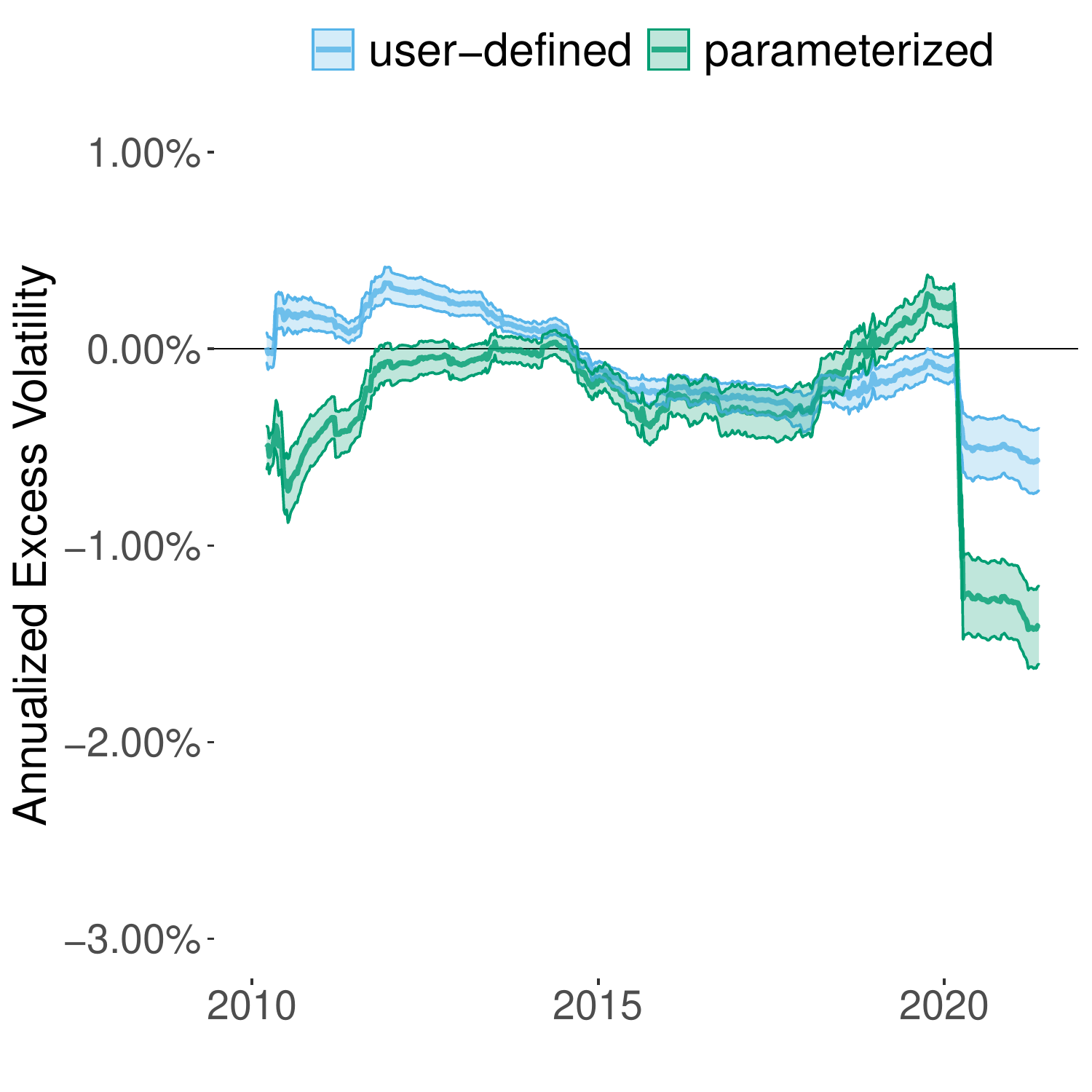}
    \caption{L2.}
  \end{subfigure}
  \begin{subfigure}[b]{0.45\linewidth}
    \includegraphics[width=\linewidth , trim={0mm 0cm 0cm 0cm},clip]{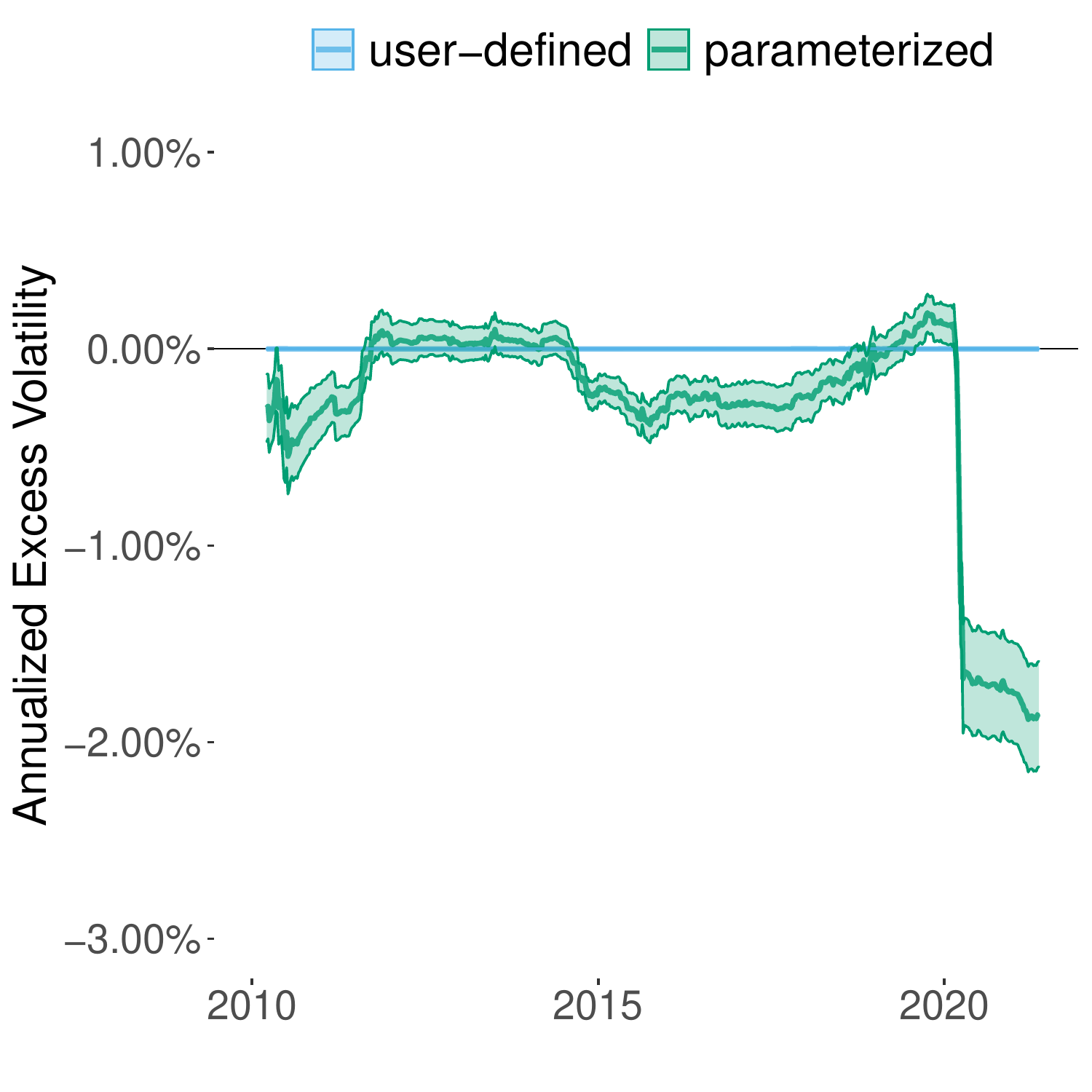}
    \caption{L2-COV.}
  \end{subfigure}
    \begin{subfigure}[b]{0.45\linewidth}
    \includegraphics[width=\linewidth , trim={0mm 0cm 0cm 0cm},clip]{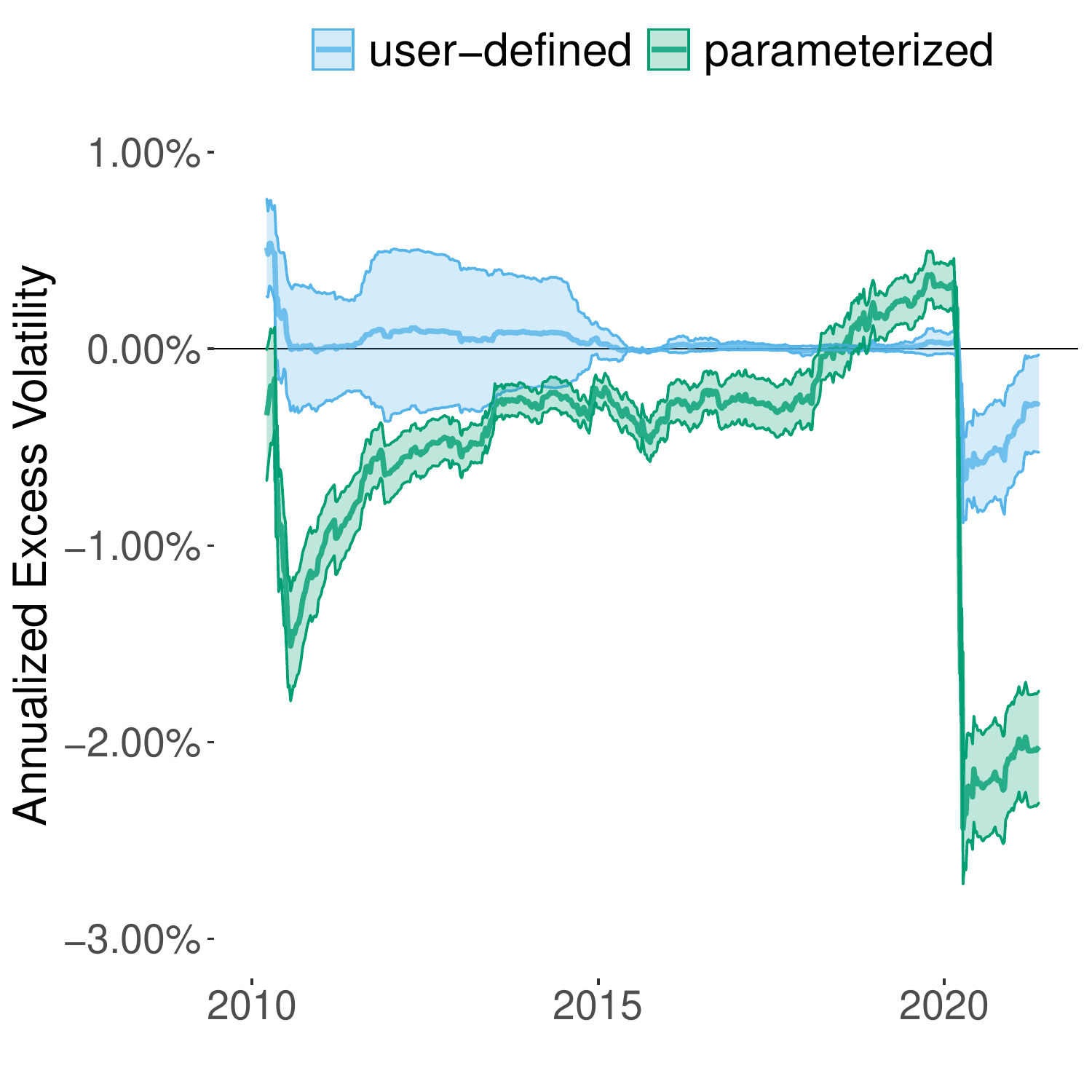}
    \caption{L1.}
  \end{subfigure}
  \begin{subfigure}[b]{0.45\linewidth}
    \includegraphics[width=\linewidth , trim={0mm 0cm 0cm 0cm},clip]{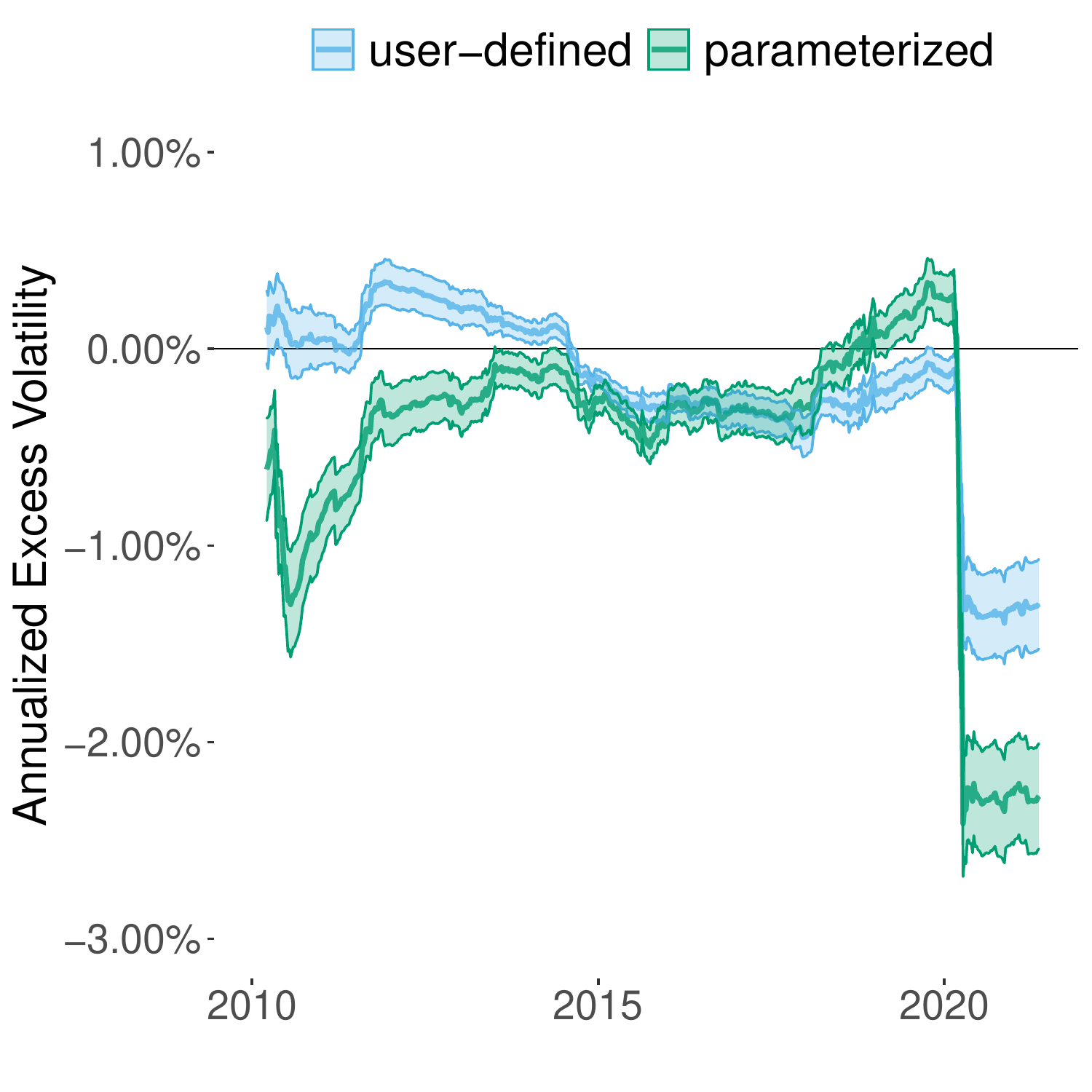}
    \caption{EN.}
  \end{subfigure}
  \caption{Out-of-sample rolling average $3$-year excess annualized volatility with $95\%$ confidence interval evaluated over $30$ randomized trials.}
  \label{fig:mv_roll}
\end{figure}

Figure \ref{fig:mv_roll} reports the out-of-sample rolling average $3$-year annualized volatility difference between the norm-penalized minimum-variance portfolios and the nominal minimum-variance portfolio. The shaded regions denote the  $95\%$ confidence interval evaluated over $30$ randomized trials. Note that values less than $0\%$ indicate that the norm-penalized portfolios produce a lower realized volatility over that time period. For user-defined models we generally observe insignificant difference in out-of-sample volatility with values oscillating around the zero line for the majority of the out-of-sample period. Notably, the L2-COV model results in a $0\%$ reduction in volatility, which is not surprising given that the penalty structure is closely related to the  covariance estimate itself. In contrast, we observe that in normal market conditions, the parameterized models result in a consistent reduction in out-of-sample volatility, with the exception of $2018-2019$ where the norm-penalized models produce marginally higher volatility. Both the user-defined and parameterized models exhibit a large reduction in realized volatility during the COVID-19 crisis of $2020-2021$. During that time the user-defined models result in a $0\% - 1.32\%$ average reduction in annualized volatility whereas the parameterized models result in a $1.27\% - 2.27\%$ annualized volatility reduction in comparison to the nominal portfolio.  

Finally in Table \ref{tab:mv}  we report the average out-of-sample economic performance for the nominal and norm-penalized minimum-variance portfolios. We also report standard errors evaluated over the $30$ randomized trials. Observe that, with the exception of the  L2-COV model, all norm-penalized models produce a lower out-of-sample average volatility; which is encouraging as this is exactly the objective for which the norm-penalty parameters are optimized. With the exception of the L2 model, all norm-penalized models produce lower average returns. The resulting impact on Sharpe ratios is mixed with the L2 and L2-COV-P models producing marginally higher Sharpe ratios and all other models producing lower out-of-sample risk-adjusted returns.

\begin{table}[h]
\centering
\begin{tabular}{l l l l l l l}
\toprule
  & Mean &$\sigma_{\text{Mean}}$& Volatility    & $\sigma_{\text{Vol}}$ & Sharpe   Ratio & $\sigma_{\text{Sharpe}}$\\
\hline
Nominal & 0.1265 & 0.0020 & 0.1411 & 0.0016 & 0.8682 & 0.0138\\
\\
L2 & 0.1267 & 0.0017 & 0.1387 & 0.0014 & 0.8852 & 0.0138\\
L2-P & 0.1236 & 0.0019 & 0.1350 & 0.0014 & 0.8862 & 0.0141\\
\\
L2-COV & 0.1265 & 0.0020 & 0.1411 & 0.0016 & 0.8681 & 0.0138\\
L2-COV-P & 0.1246 & 0.0021 & 0.1334 & 0.0013 & 0.9039 & 0.0163\\
\\
L1 & 0.1151 & 0.0026 & 0.1403 & 0.0018 & 0.7906 & 0.0138\\
L1-P & 0.1048 & 0.0020 & 0.1320 & 0.0011 & 0.7626 & 0.0143\\
\\
EN & 0.1163 & 0.0016 & 0.1357 & 0.0012 & 0.8271 & 0.0123\\
EN-P & 0.1081 & 0.0020 & 0.1313 & 0.0012 & 0.7923 & 0.0150\\
\bottomrule
\end{tabular}
\caption{Average and standard errors  of out-of-sample economic performance for the nominal and norm-penalized minimum-variance portfolios evaluated over $30$ randomized trials.  }
\label{tab:mv}
\end{table}

\subsection{Experiment 2: global futures data}\label{sec:exp_2}
We consider an unconstrained mean-variance optimization $(\mathbb{S} = \mathbb{R}^{d_z})$  on global futures data, described in Table \ref{table:futures}. The norm-penalty parameters are determined by solving Program $\eqref{eq:stoch_discrete}$  to minimize the  Sharpe ratio decision cost, $c_{\text{SR}}$. Asset mean returns are estimated by least-squares according to Equation $\eqref{eq:y_hat}$. The feature variables, $\toi{\bx} \in \mathbb{R}^5$, are the average trailing $1$-month risk-adjusted return for each market, averaged by asset class. The asset covariance matrix, $\toi{\hat{ \bV } } $, is estimated by a simple moving average covariance evaluated over the trailing $52$ weeks and the risk-aversion parameter is set to $\delta = 20$. In each trial we randomly draw $n=20$ assets from the asset universe. 


\begin{figure}[]
  \centering
    \includegraphics[width=0.75\linewidth , trim={0mm 0cm 0cm 0cm},clip]{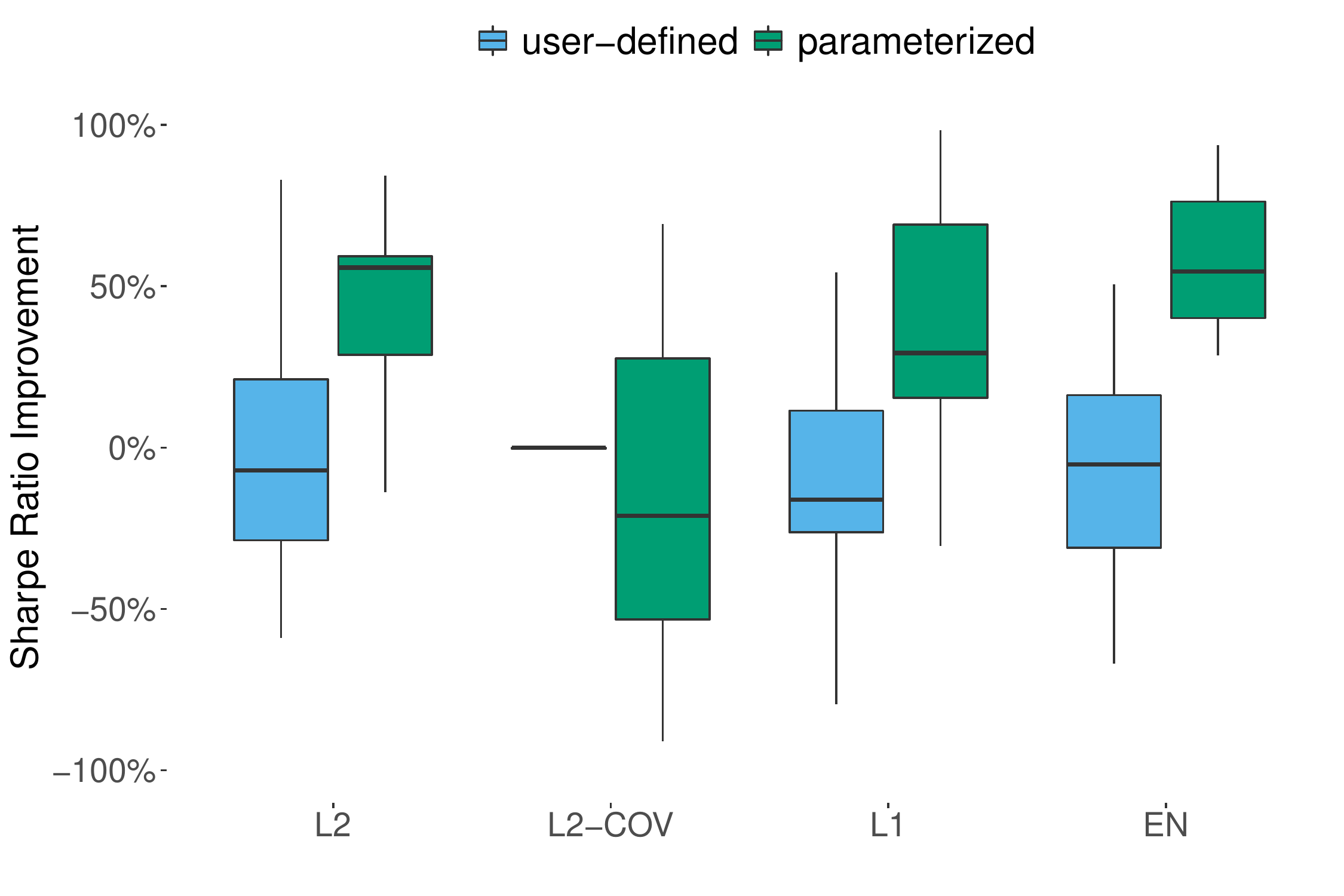}
   \caption{Out-of-sample Percent Sharpe ratio improvement of norm-penalized mean-variance portfolios evaluated over $30$ randomized trials}
  \label{fig:msr_box}
\end{figure}

Figure \ref{fig:msr_box} reports the out-of-sample percent improvement in Sharpe ratios of the norm-penalized MVO portfolios relative to the nominal MVO portfolio. Note here that values greater than $0\%$ indicate an increase in Sharpe ratio, which we seek to maximize. With the exception of the L2-COV-P model, all other norm-penalized portfolios increase the average out-of-sample Sharpe ratio with median percent improvements ranging from $0\% - 45.7\%$.  We note that the improvement in the out-of-sample portfolio objective is not as consistent as the reduction in variance exhibited by the stocks data experiment, with many trials resulting in lower out-of-sample Sharpe ratios. As before we note that the parameterized models result in a larger improvement in out-of-sample portfolio objectives in comparison to their user-defined counterparts. The parameterized models produce median out-of-sample Sharpe ratio improvements ranging from $-8.1\% - 45.7\%$ whereas the user-defined models produce median out-of-sample Sharpe ratio improvements ranging from $0\% - 8.5\%$. Lastly, parameterized models exhibit a larger variance in out-of-sample Sharpe ratio, suggesting that the models may be overfitting the in-sample data. 

\begin{figure}[]
  \centering
  \begin{subfigure}[b]{0.45\linewidth}
    \includegraphics[width=\linewidth , trim={0mm 0cm 0cm 0cm},clip]{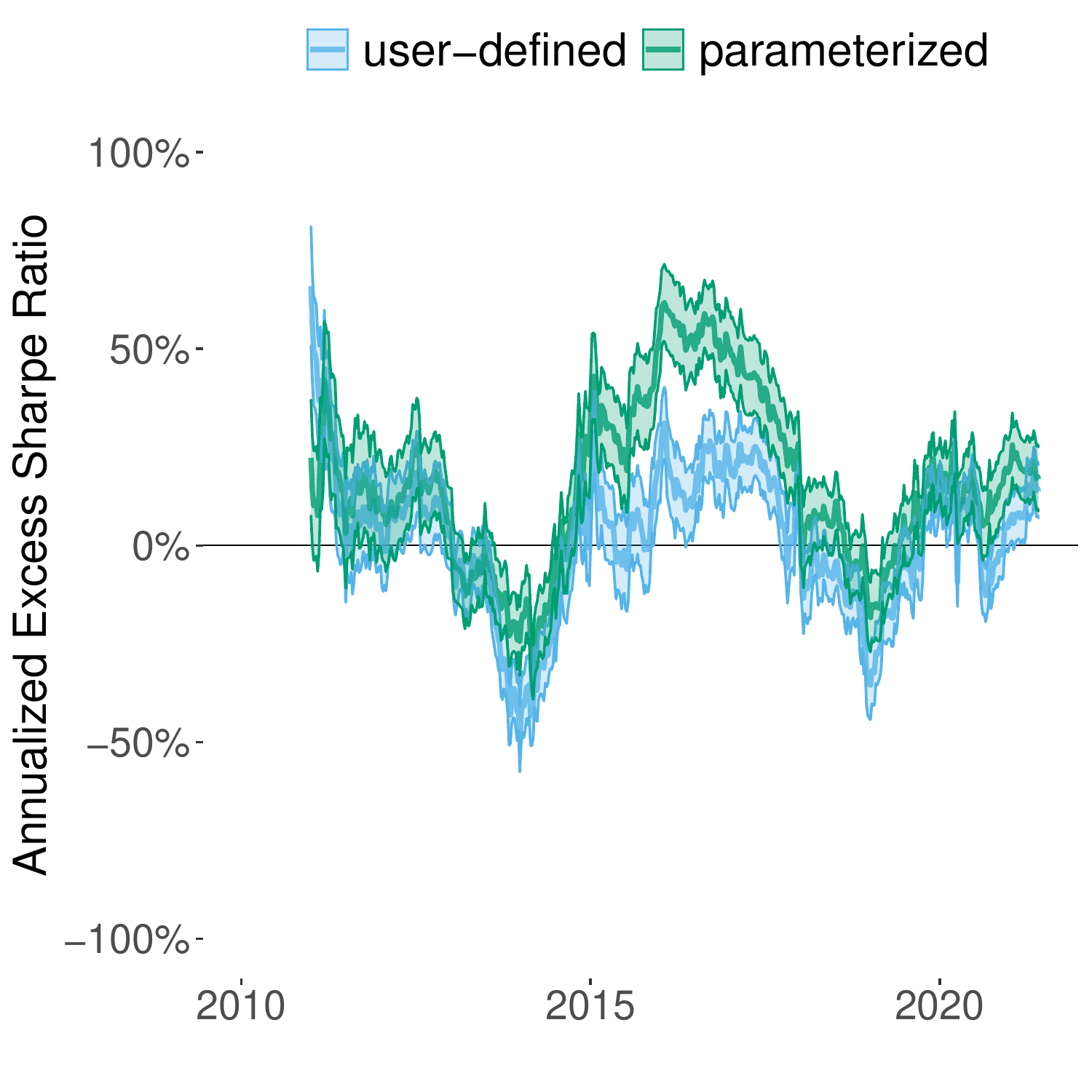}
    \caption{L2.}
  \end{subfigure}
  \begin{subfigure}[b]{0.45\linewidth}
    \includegraphics[width=\linewidth , trim={0mm 0cm 0cm 0cm},clip]{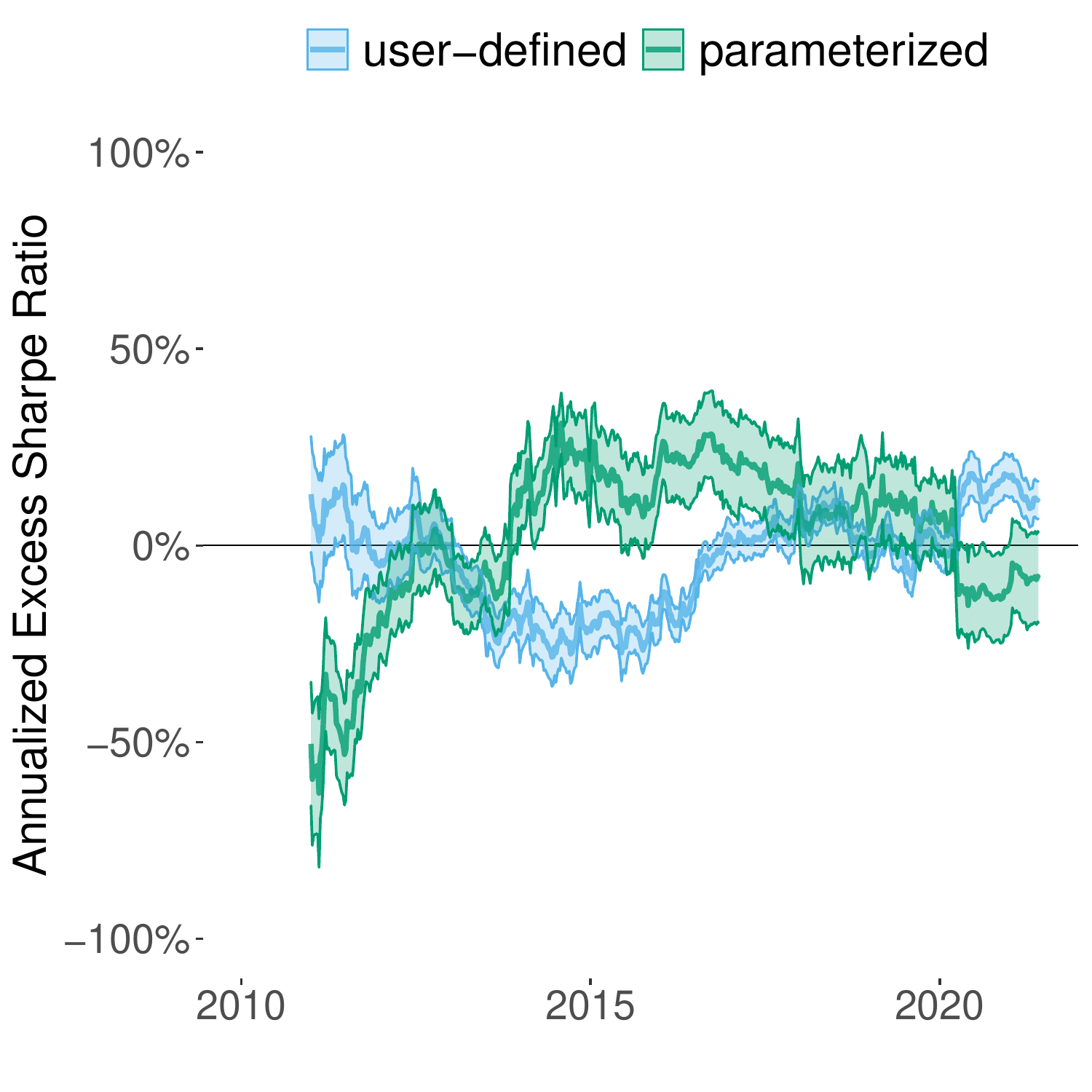}
    \caption{L2-COV.}
  \end{subfigure}
    \begin{subfigure}[b]{0.45\linewidth}
    \includegraphics[width=\linewidth , trim={0mm 0cm 0cm 0cm},clip]{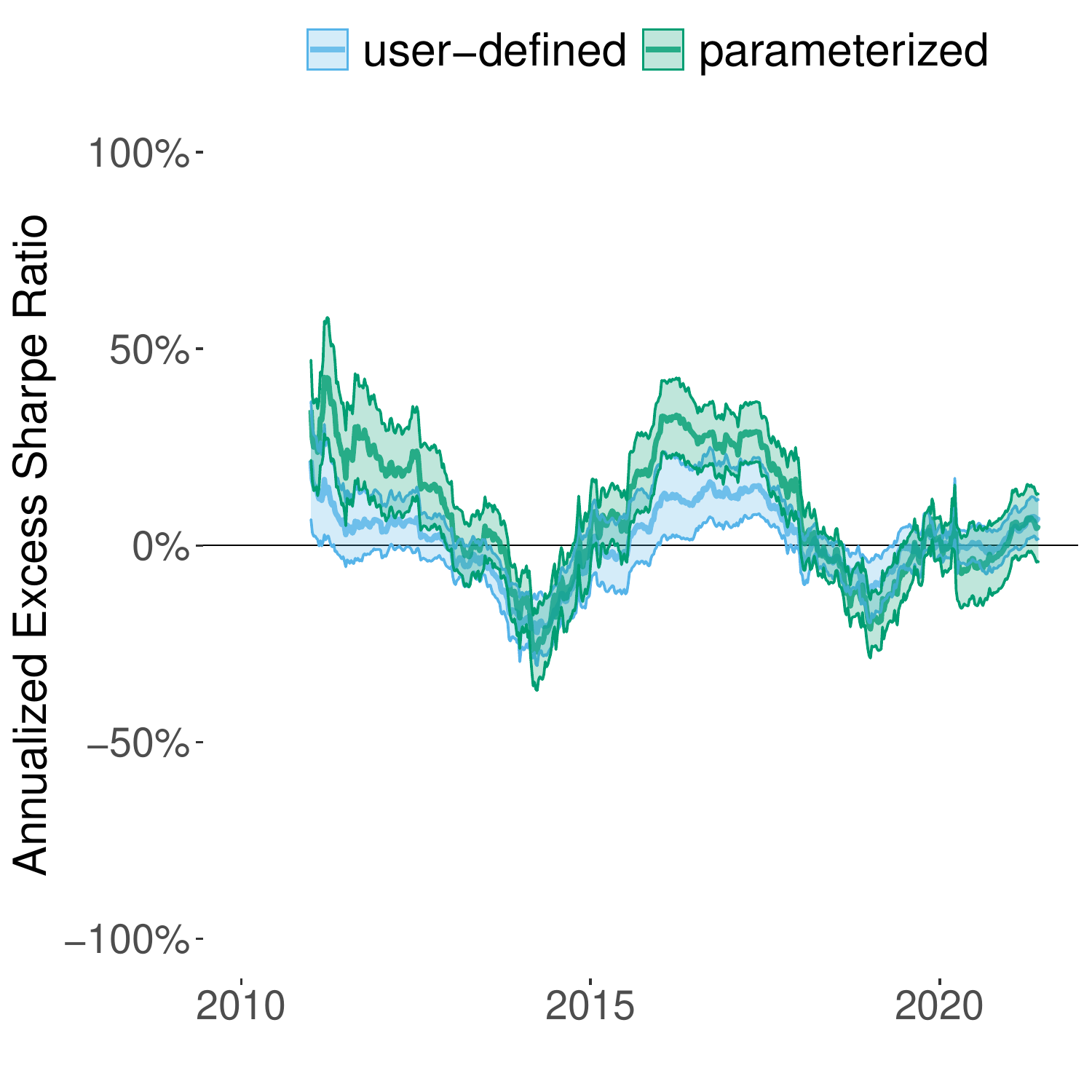}
    \caption{L1.}
  \end{subfigure}
  \begin{subfigure}[b]{0.45\linewidth}
    \includegraphics[width=\linewidth , trim={0mm 0cm 0cm 0cm},clip]{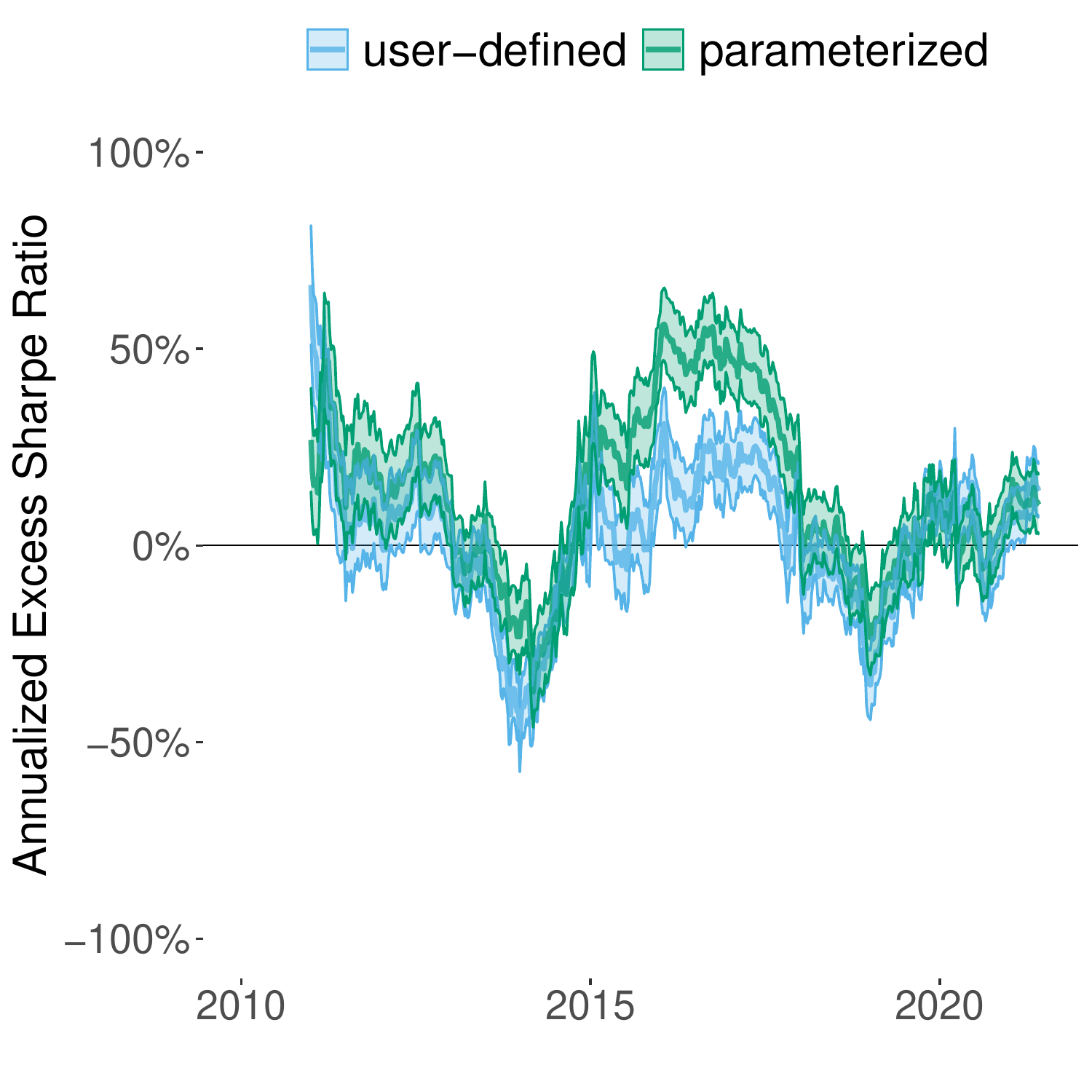}
    \caption{EN.}
  \end{subfigure}
  \caption{Out-of-sample rolling average $3$-year excess annualized Sharpe ratio with $95\%$ confidence interval evaluated over $30$ randomized trials. .}
  \label{fig:msr_roll}
\end{figure}

Figure \ref{fig:msr_roll} reports the out-of-sample rolling average $3$-year annualized Sharpe ratio difference and $95\%$ confidence interval between the norm-penalized MVO portfolios and the nominal MVO portfolio. Note that values greater than $0\%$ indicate that the norm-penalized portfolios produce a higher realized Sharpe ratio over that time period. As expected,  the parameterized models produce higher and more consistent improvements in Sharpe ratio in comparison to the corresponding user-defined models. The difference, however, is not statistically significant and displays considerable overlap in confidence intervals. Moreover, during both $2013-2014$ and $2018-2019$ we observe that the parameterized norm-penalized MVO portfolios produce lower out-of-sample Sharpe ratios in comparison to the nominal MVO.

Lastly in Table \ref{tab:msr}  we report the average out-of-sample economic performance and corresponding standard errors for the nominal and norm-penalized MVO portfolios.   As expected, the norm-penalized portfolios exhibit a reduction in portfolio volatility. With the exception of the L2-COV-P model, the reduction in volatility exceeds the corresponding reduction in mean return, thus producing higher average out-of-sample Sharpe ratios. The average Sharpe ratio improvement ranges from $4.8\%-48.3\%$. Again these results are encouraging as the norm-penalty parameters are directly optimized to maximize the Sharpe ratio. In general we find that the L2-P and EN-P models produce significant  average Sharpe ratio improvements, with average values of $0.4326 - 0.4658$ in comparison to the nominal Sharpe ratio of $0.3184$.
\begin{table}[h]
\centering
\begin{tabular}{l l l l l l l}
\toprule
  & Mean &$\sigma_{\text{Mean}}$& Volatility    & $\sigma_{\text{Vol}}$ & Sharpe   Ratio & $\sigma_{\text{Sharpe}}$\\
\hline
Nominal & 0.2118 & 0.0236 & 0.6393 & 0.0382 & 0.3184 & 0.0314\\
\\
L2 & 0.0581 & 0.0047 & 0.1652 & 0.0065 & 0.3519 & 0.0227\\
L2-P & 0.0945 & 0.0061 & 0.2028 & 0.0056 & 0.4658 & 0.0267\\
\\
L2-COV & 0.1092 & 0.0120 & 0.3385 & 0.0255 & 0.3296 & 0.0315\\
L2-COV-P & 0.1090 & 0.0155 & 0.3469 & 0.0117 & 0.3108 & 0.0419\\
\\
L1 & 0.1733 & 0.0212 & 0.4954 & 0.0436 & 0.3425 & 0.0273\\
L1-P & 0.1627 & 0.0190 & 0.4243 & 0.0357 & 0.3757 & 0.0346\\
\\
EN & 0.0584 & 0.0048 & 0.1655 & 0.0067 & 0.3530 & 0.0227\\
EN-P & 0.0857 & 0.0064 & 0.1983 & 0.0068 & 0.4326 & 0.0289\\
\bottomrule
\end{tabular}
\caption{Average and standard errors  of out-of-sample economic performance for the nominal and norm-penalized MVO portfolios evaluated over $30$ randomized trials.  }
\label{tab:msr}
\end{table}

\section{Conclusion and future work}\label{sec:conclusion}
In this paper we augmented the nominal mean-variance program with a convex combination of parameterized $L_1$ and $L_2$-norm penalty functions. A traditional approach would pre-specify the penalty structure and then optimize for the amount of regularization. Alternatively, we presented a data-driven framework for optimizing the parameterized penalty structures such that the resulting penalized program induces optimal portfolio decision-making. The parameter estimation problem is cast as a bi-level program and optimized by first-order gradient descent. Specifically, we structured the learning process an an end-to-end neural network with differentiable quadratic programming layers. We proposed a novel technique for differentiating the optimal solution, $\bz^*$, with respect to the parameterized $L_1$-norm, which solves the dual program - a box-constrained QP - and then recovers the primal solution from the optimal dual variables. 

We presented several experiments using both US stocks and global futures data and demonstrated the benefit of the data-driven approach. In general we find that the parameterized norm-penalty models, in which the exact penalty structure is learned from the data, provide improved out-of-sample portfolio objectives in comparison to the nominal and user-defined norm-penalized portfolios. Our experiments should be interpreted as a proof-of-concept and we acknowledge that further testing with alternative datasets,  and under varying prediction model and portfolio constraint assumptions is required in order to better determine the efficacy of the data-driven solution. In particular, we note that the parameterized models, if poorly specified, risk overfitting the training dataset.  We believe that a stochastic optimization framework whereby model parameters are optimized in order to minimize a cross-validation cost  may attenuate the propensity for model overfit and is an active area of research. 

\bibliographystyle{plainnat}
\bibliography{Bibliography/Bibliography}

\begin{thebibliography}{51}
\providecommand{\natexlab}[1]{#1}
\providecommand{\url}[1]{\texttt{#1}}
\expandafter\ifx\csname urlstyle\endcsname\relax
  \providecommand{\doi}[1]{doi: #1}\else
  \providecommand{\doi}{doi: \begingroup \urlstyle{rm}\Url}\fi

\bibitem[Agrawal et~al.(2019)Agrawal, Amos, Barratt, Boyd, Diamond, and
  Kolter]{Agrawal2019}
Akshay Agrawal, Brandon Amos, Shane Barratt, Stephen Boyd, Steven Diamond, and
  J.~Zico Kolter.
\newblock Differentiable convex optimization layers.
\newblock In \emph{Advances in Neural Information Processing Systems},
  volume~32, pages 9562--9574. Curran Associates, Inc., 2019.

\bibitem[Amos and Kolter(2017)]{Amos2017}
Brandon Amos and J.~Zico Kolter.
\newblock Optnet: Differentiable optimization as a layer in neural networks,
  2017.
\newblock URL \url{https://arxiv.org/abs/1703.00443}.

\bibitem[Amos et~al.(2019)Amos, Rodriguez, Sacks, Boots, and Kolter]{Amos2019}
Brandon Amos, Ivan Dario~Jimenez Rodriguez, Jacob Sacks, Byron Boots, and
  J.~Zico Kolter.
\newblock Differentiable mpc for end-to-end planning and control, 2019.
\newblock URL \url{http://arxiv.org/abs/1810.13400}.

\bibitem[Ben-Tal et~al.(2009)Ben-Tal, Ghaoui, and ANemirovski]{Ben2009}
A~Ben-Tal, L~El Ghaoui, and ANemirovski.
\newblock \emph{Robust Optimization}.
\newblock Princeton University Press, 2009.

\bibitem[Ben-Tal and Nemirovski(2000)]{Ben2000}
Aharon Ben-Tal and Arkadi Nemirovski.
\newblock Robust solutions of linear programming problems contaminated with
  uncertain data.
\newblock \emph{Mathematical Programming}, 88:\penalty0 411--424, 01 2000.
\newblock \doi{10.1007/PL00011380}.

\bibitem[Bertsimas and Kallus(2020)]{Bert2020}
Dimitris Bertsimas and Nathan Kallus.
\newblock From predictive to prescriptive analytics.
\newblock \emph{Management Science}, 66\penalty0 (3):\penalty0 1025--1044,
  2020.

\bibitem[Bertsimas and Sim(2004)]{Bert2004}
Dimitris Bertsimas and Melvyn Sim.
\newblock The price of robustness.
\newblock \emph{Operations Research}, 52:\penalty0 35--53, 02 2004.
\newblock \doi{10.1287/opre.1030.0065}.

\bibitem[Bertsimas et~al.(2014)Bertsimas, Gupta, and Kallus]{Bert2014}
Dimitris Bertsimas, Vishal Gupta, and Nathan Kallus.
\newblock Data-driven robust optimization, 2014.

\bibitem[Black and Litterman(1991)]{Black1991}
F.~Black and R.~Litterman.
\newblock Asset allocation combining investor views with market equilibrium.
\newblock \emph{Journal of Fixed Income}, 1\penalty0 (2):\penalty0 7--18, 1991.

\bibitem[Boyd and Vandenberghe(2004)]{Boyd2004}
Stephen Boyd and Lieven Vandenberghe.
\newblock \emph{Convex Optimization}.
\newblock Cambridge University Press, 2004.
\newblock \doi{10.1017/CBO9780511804441}.

\bibitem[Butler and Kwon(2021{\natexlab{a}})]{Butler2021ADMM}
Andrew Butler and Roy Kwon.
\newblock Efficient differentiable quadratic programming layers: an admm
  approach, 2021{\natexlab{a}}.
\newblock URL \url{https://arxiv.org/abs/2112.07464}.

\bibitem[Butler and Kwon(2021{\natexlab{b}})]{Butler2021IPOb}
Andrew Butler and Roy Kwon.
\newblock Covariance estimation for risk-based portfolio optimization: an
  integrated approach.
\newblock \emph{Journal of Risk}, 24\penalty0 (2), December 2021{\natexlab{b}}.

\bibitem[Butler and Kwon(2021{\natexlab{c}})]{Butler2020IPO}
Andrew Butler and Roy~H. Kwon.
\newblock Integrating prediction in mean-variance portfolio optimization,
  2021{\natexlab{c}}.
\newblock URL \url{https://arxiv.org/abs/2102.09287}.

\bibitem[Chevalier et~al.(2022)Chevalier, Coqueret, and Raffinot]{Chev2022}
Guillaume Chevalier, Guillaume Coqueret, and Thomas Raffinot.
\newblock Supervised portfolios.
\newblock \emph{SSRN Electronic Journal}, 4 2022.
\newblock URL \url{https://ssrn.com/abstract=3954109}.

\bibitem[Chopra and Ziemba(1993)]{Chopra1993}
Vijay~Kumar. Chopra and William~T. Ziemba.
\newblock The effect of errors in means, variances, and covariances on optimal
  portfolio choice.
\newblock \emph{The Journal of Portfolio Management}, 19\penalty0 (2):\penalty0
  6--11, 1993.
\newblock ISSN 0095-4918.
\newblock \doi{10.3905/jpm.1993.409440}.
\newblock URL \url{https://jpm.pm-research.com/content/19/2/6}.

\bibitem[DeMiguel et~al.(2009)DeMiguel, Garlappi, Uppal, and
  Nogales]{Demiguel2009}
Victor DeMiguel, Lorenzo Garlappi, Raman Uppal, and L~Nogales.
\newblock A generalized approach to portfolio optimization: improving
  performance by constraining portfolio norms.
\newblock \emph{Management Science}, 55\penalty0 (5):\penalty0 798--812, 2009.

\bibitem[Donti et~al.(2017)Donti, Amos, and Kolter]{Donti2017}
Priya Donti, Brandon Amos, and J.~Zico Kolter.
\newblock Task-based end-to-end model learning in stochastic optimization.
\newblock In I.~Guyon, U.~V. Luxburg, S.~Bengio, H.~Wallach, R.~Fergus,
  S.~Vishwanathan, and R.~Garnett, editors, \emph{Advances in Neural
  Information Processing Systems}, volume~30, pages 5484 -- 5494. Curran
  Associates, Inc., 2017.

\bibitem[Drees and Starica(2002)]{Drees2002}
Holger Drees and Catalin Starica.
\newblock A simple non-stationary model for stock returns, 2002.

\bibitem[Elmachtoub and Grigas(2017)]{Elma2020}
Adam Elmachtoub and Paul Grigas.
\newblock Smart `predict, then optimize'.
\newblock \emph{Management Science}, 10 2017.
\newblock \doi{10.1287/mnsc.2020.3922}.

\bibitem[Elmachtoub et~al.(2020)Elmachtoub, Liang, and McNellis]{Elma2020b}
Adam~N. Elmachtoub, Jason Cheuk~Nam Liang, and Ryan McNellis.
\newblock Decision trees for decision-making under the predict-then-optimize
  framework, 2020.
\newblock URL \url{http://arxiv.org/abs/2003.00360}.

\bibitem[Engle(1982)]{Engle1982}
Robert~F. Engle.
\newblock Autoregressive conditional heteroskedasticity with estimates of the
  variance of uk inflation.
\newblock \emph{Econometrica}, 50\penalty0 (1):\penalty0 987--1008, 1982.

\bibitem[Feng and Simon(2017)]{Feng2017}
Jean Feng and Noah Simon.
\newblock Gradient-based regularization parameter selection for problems with
  non-smooth penalty functions, 2017.

\bibitem[Goldfarb and Iyengar(2003)]{Goldfarb2003}
D.~Goldfarb and G.~Iyengar.
\newblock Robust portfolio selection problems.
\newblock \emph{Mathematics of Operations Research}, 28\penalty0 (1):\penalty0
  1--38, 2003.

\bibitem[Grigas et~al.(2021)Grigas, Qi, Zuo-Jun, and Shen]{Grigas2021}
Paul Grigas, Meng Qi, Zuo-Jun, and Shen.
\newblock Integrated conditional estimation-optimization, 2021.
\newblock URL \url{http://arxiv.org/abs/2110.12351}.

\bibitem[Hastie et~al.(2001)Hastie, Tibshirani, and Friedman]{Hastie2001}
Trevor Hastie, Robert Tibshirani, and Jerome Friedman.
\newblock \emph{The Elements of Statistical Learning}.
\newblock Springer Series in Statistics. Springer New York Inc., New York, NY,
  USA, 2001.

\bibitem[Ho et~al.(2015)Ho, Sun, and Xin]{Ho2015}
Michael Ho, Zheng Sun, and Jack Xin.
\newblock Weighted elastic net penalized mean-variance portfolio design and
  computation.
\newblock \emph{SIAM Journal on Financial Mathematics}, 6\penalty0
  (1):\penalty0 1220--1244, 2015.

\bibitem[Ichnowski et~al.(2021)Ichnowski, Jain, Stellato, Banjac, Luo,
  Borrelli, Gonzalez, Stoica, and Goldberg]{Ich2021}
Jeffrey Ichnowski, Paras Jain, Bartolomeo Stellato, Goran Banjac, Michael Luo,
  Francesco Borrelli, Joseph~E. Gonzalez, Ion Stoica, and Ken Goldberg.
\newblock Accelerating quadratic optimization with reinforcement learning,
  2021.
\newblock URL \url{https://arxiv.org/abs/2107.10847}.

\bibitem[Jobson and Korkie(1980)]{Jobson1980}
J.D. Jobson and Bob Korkie.
\newblock Estimation for markowitz efficient portfolios.
\newblock \emph{Journal of the American Statistical Association}, pages 544 --
  555, 1980.

\bibitem[Jobson and Korkie(1982)]{Jobson1982}
J.D. Jobson and Bob Korkie.
\newblock Potential performance and tests of portfolio efficiency.
\newblock \emph{Journal of Financial Economics}, 10\penalty0 (4):\penalty0 433
  -- 466, 1982.

\bibitem[Kim et~al.(2008)Kim, Koh, Lustig, Boyd, and Gorinevsky]{Kim2008}
Seung-Jean Kim, K.~Koh, M.~Lustig, Stephen Boyd, and Dimitry Gorinevsky.
\newblock An interior-point method for large-scale l1-regularized least
  squares.
\newblock \emph{Selected Topics in Signal Processing, IEEE Journal of},
  1:\penalty0 606 -- 617, 01 2008.
\newblock \doi{10.1109/JSTSP.2007.910971}.

\bibitem[Kingma and Ba(2014)]{Kingma2015}
Diederik~P. Kingma and Jimmy Ba.
\newblock Adam: A method for stochastic optimization, 2014.
\newblock URL \url{https://arxiv.org/abs/1412.6980}.

\bibitem[Ledoit and Wolf(2004{\natexlab{a}})]{Ledoit2003}
Olivier Ledoit and Michael Wolf.
\newblock Honey, i shrunk the sample covariance matrix.
\newblock \emph{The Journal of Portfolio Management}, 30\penalty0 (4):\penalty0
  110--119, 2004{\natexlab{a}}.

\bibitem[Ledoit and Wolf(2004{\natexlab{b}})]{Ledoit2004}
Olivier Ledoit and Michael Wolf.
\newblock A well-conditioned estimator for large-dimensional covariance
  matrices.
\newblock \emph{Journal of Multivariate Analysis}, 88\penalty0 (2):\penalty0
  365--411, 2004{\natexlab{b}}.

\bibitem[Ledoit and Wolf(2012)]{Ledoit2012}
Olivier Ledoit and Michael Wolf.
\newblock Nonlinear shrinkage estimation of large-dimensional covariance
  matrices.
\newblock \emph{The Annals of Statistics}, 40\penalty0 (2):\penalty0
  1024--1060, 2012.

\bibitem[Lorraine and Duvenaud(2018)]{Lorraine2018}
Jonathan Lorraine and David Duvenaud.
\newblock Stochastic hyperparameter optimization through hypernetworks, 2018.

\bibitem[Mandi and Guns(2020)]{Mandi2020}
Jayanta Mandi and Tias Guns.
\newblock Interior point solving for lp-based prediction+optimisation, 2020.
\newblock URL \url{http://arxiv.org/abs/2010.13943}.

\bibitem[Mandi et~al.(2019)Mandi, Demirovic, Stuckey, and Guns]{Mandi2019}
Jaynta Mandi, Emir Demirovic, Peter.~J Stuckey, and Tias Guns.
\newblock Smart predict-and-optimize for hard combinatorial optimization
  problems, 2019.
\newblock URL \url{http://arxiv.org/abs/1911.10092}.

\bibitem[Markowitz(1952)]{Markowitz1952}
H.~Markowitz.
\newblock Portfolio selection.
\newblock \emph{Journal of Finance}, 7\penalty0 (1):\penalty0 77--91, 1952.

\bibitem[Michaud and Michaud(2008{\natexlab{a}})]{Michaud2008a}
Richard Michaud and Robert Michaud.
\newblock Efficient asset management: A practical guide to stock portfolio
  optimization and asset allocation.
\newblock \emph{New York: Oxford University Press}, 1, 2008{\natexlab{a}}.

\bibitem[Michaud and Michaud(2008{\natexlab{b}})]{Michaud2008b}
Richard Michaud and Robert Michaud.
\newblock Estimation error and portfolio optimization: A resampling solution.
\newblock \emph{Journal of Investment Management}, 6\penalty0 (1):\penalty0
  8--28, 2008{\natexlab{b}}.

\bibitem[Pedregosa(2016)]{Ped2016}
Fabian Pedregosa.
\newblock Hyperparameter optimization with approximate gradient, 2016.

\bibitem[Rumelhart et~al.(1986)Rumelhart, Hinton, and Williams]{Hinton1986}
David~E. Rumelhart, Geoffrey~E. Hinton, and Ronald~J. Williams.
\newblock Learning representations by back-propagating errors.
\newblock \emph{Nature}, 323:\penalty0 533--536, 1986.

\bibitem[Starica and Granger(2005)]{Starica2005}
Catalin Starica and Clive Granger.
\newblock Nonstationarities in stock returns.
\newblock \emph{The Review of Economics and Statistics}, 87\penalty0
  (3):\penalty0 503--522, 2005.

\bibitem[Stein(1956)]{Stein1956}
Charles Stein.
\newblock Inadmissibility of the usual estimator for the mean of a multivariate
  normal distribution.
\newblock In \emph{Proceedings of the Third Berkeley Symposium on Mathematical
  Statistics and Probability, Volume 1: Contributions to the Theory of
  Statistics}, pages 197--206, Berkeley, Calif., 1956. University of California
  Press.
\newblock URL \url{https://projecteuclid.org/euclid.bsmsp/1200501656}.

\bibitem[Tibshirani(1996)]{Tibshirani1996}
Robert Tibshirani.
\newblock Regression shrinkage and selection via the lasso.
\newblock \emph{Journal of the Royal Statistical Society}, 58\penalty0 (1),
  1996.

\bibitem[Tikhonov(1963)]{Tikhonov1963}
A.~N. Tikhonov.
\newblock Solution of incorrectly formulated problemsand the regularization
  method.
\newblock \emph{Soviet Mathematics}, pages 1035--1038, 1963.

\bibitem[Tutuncu and Koenig(2004)]{Tu2004}
R.~Tutuncu and M.~Koenig.
\newblock Robust asset allocation.
\newblock \emph{Annals of Operations Research}, 132:\penalty0 157--187, 2004.

\bibitem[Uysal et~al.(2021)Uysal, Li, and Mulvey]{Uysal2021}
Ayse~Sinem Uysal, Xiaoyue Li, and John~M. Mulvey.
\newblock End-to-end risk budgeting portfolio optimization with neural
  networks, 2021.
\newblock URL \url{http://arxiv.org/abs/2107.04636}.

\bibitem[Yin et~al.(2019)Yin, Perchet, and Soupe]{Yin2019}
Chenyang Yin, Romain Perchet, and Francois Soupe.
\newblock A practical guide to robust portfolio optimization.
\newblock \emph{SSRN Electronic Journal}, 1 2019.
\newblock \doi{10.2139/ssrn.3490680}.

\bibitem[Zhu et~al.(2009)Zhu, Coleman, and li]{Zhu2009}
Lei Zhu, Thomas Coleman, and Yuying li.
\newblock Min-max robust and cvar robust mean-variance portfolios.
\newblock \emph{Journal of Risk}, 11, 03 2009.
\newblock \doi{10.21314/JOR.2009.191}.

\bibitem[Zou and Hastie(2005)]{Zou2005}
Hui Zou and Trevor Hastie.
\newblock Regularization and variable selection via the elastic net.
\newblock \emph{Journal of the Royal Statistical Society. Series B (Statistical
  Methodology)}, 67\penalty0 (2):\penalty0 301--320, 2005.

\end{thebibliography}

\appendix

\section{Implementation details } \label{sec:app_imp_details}

The norm-penalty coefficients, $\bgamma_1$ and $\bgamma_2$ are constrained to the positive real orthant which is satisfied by applying following exponential transformation:
$$\bgamma_1 \leftarrow e^{\bgamma_1} \qquad \text{and} \bgamma_2 \leftarrow e^{\bgamma_1}.$$
Coefficients are initialized at $\bgamma_1 = -4$ and $\bgamma_2 =-4$. The diagonal matrices  ${\bE(\btheta_1)= \diag(\btheta_1)}$ and ${\bD(\btheta_2) = \diag(\btheta_2)}$ are constrained to be positive semidefinite, which we implement by applying a relu transformation:
$$  \btheta_1 \leftarrow \max( \btheta_1,0) \quad \text{and} \quad \btheta_2   \leftarrow  \max( \btheta_2,0)  $$
The parameters are initialized randomly from the standard uniform distributions:
$$
\btheta_1 \in [0,1]^{d_z}, \quad \btheta_2 \in [0,1]^{d_z}
$$
All parameters are optimized by applying the ADAM gradient descent routine \citep{Kingma2015}.  We find that using $100\%$ of the training observations with a learning rate of $0.10$ results in numerically stable and consistent parameter estimation. All experiments were conducted on an Apple Mac Pro computer ($2.7$ GHz $12$-Core Intel Xeon E5,$128$ GB $1066$ MHz DDR3 RAM) running macOS `Catalina'. The software was written using the R programming language (version $4.0.0$) and  torch (version $0.6.0$).

The code snippet, presented  below demonstrates how to implement the optimization of the norm-penalty parameters as an end-to-end trainable neural network for the example provided in Section \ref{sec:network} . We refer to the open-source R package \textit{lqp} (learning quadratic programs), available here: \textit{https://github.com/butl3ra/lqp}, for more details. Note that the norm-penalty contains  a parameterized $L_1$-norm and therefore we implement the dual solver with fixed-point implicit differentiation, as described in Section \ref{sec:dual}.
\newpage
\begin{python}
# --- setup solver control:
control = nn_qp_control(solver = 'con_l1', backprop = 'fixed_point')
# --- p_model:
p_model = nn_constant(value = torch_zeros(c(1,n_y,1)))
# --- Q_model:
Q_model = nn_constant(value = Q_is_t)
# --- A_model and b_model:
A_model = nn_constant(value = torch_ones(c(1,1,n_y)))
b_model = nn_constant(value = torch_ones(c(1)))
# --- G_model and h_model:
G_model = nn_constant(value = -torch_eye(n_y)$unsqueeze(1))
h_model = nn_constant(value = torch_zeros(c(1,n_y,1)) )
# --- E_model and lambda_1:
E = nn_parameter(as_torch_tensor(diag(runif(n_y)),requires_grad = TRUE))
E_model = nn_sequential(nn_constant(value = E),nn_unsqueeze(dim = 1),nn_relu())
lambda_1 = nn_parameter(as_torch_tensor(-4,requires_grad = TRUE))
lambda_1_model = nn_sequential(nn_constant(value = lambda_1), nn_exp())
# --- D_model and lambda_2:
D_model = nn_quad_form_const(in_features =  n_x, out_features=n_y, x_mat = Q_x)
lambda_2 = nn_parameter(as_torch_tensor(-4,requires_grad = TRUE))
lambda_2_model = nn_sequential(nn_constant(value = lambda_2), nn_exp())

# --- norm-penalized qp model:
model = nn_lqp(Q_model = Q_model, p_model = p_model,  A_model = A_model,
			 				 b_model = b_model, G_model = G_model,  h_model = h_model,
							 D_model = D_model, lambda_2_model = lambda_2_model,
							 E_model = E_model, lambda_1_model = lambda_1_model,
							 control = control)
optimizer <- optim_adam(model$parameters, lr = 0.10)
#---------------- training loop -------------------
for (t in 1:100 ) {
  # --- forward pass
  z = model(x = x_is_t)
  z = torch_sum_1(z,dim = 2)
  # --- compute loss
  loss <-  nnf_var_loss(z = z, y_is_t$unsqueeze(3))
  # --- backpropagation:
  optimizer$zero_grad()
  loss$backward()
  # ---  update params:
  optimizer$step()
}

\end{python}

\newpage
\section{Data summary } \label{sec:app_data}
\begin{table}[h]
\centering
\scriptsize
 \begin{tabular}{  l  c c c c c c c c }
 \hline
{\bf{GICS Sector}} & & & &  {\bf{Stock Symbols}}\\
\hline
Communication &   CBB & CMCSA & DIS & FOX & IPG & LUMN & MDP & NYT \\
      &    T & VOD & VZ\\
Consumer      &    BBY & CBRL & CCL & F & GPC & GPS & GT & HAS \\
Discretionary &    HD & HOG & HRB & JWN & LB & LEG & LEN & LOW \\
              &    MCD & NKE & NVR & NWL & PHM & PVH & ROST & TGT\\
              &    TJX & VFC & WHR & WWW\\
Consumer      &   ADM & ALCO & CAG & CASY & CHD & CL & CLX & COST\\
Staples       &   CPB & FLO & GIS & HSY & K & KMB & KO & KR \\
              &   MO & PEP & PG & SYY & TAP &  TR & TSN & UVV\\
              & WBA  & WMK & WMT\\
Energy        &   AE & APA  & BKR & BP & COP & CVX & EOG & HAL\\
              &   HES & MRO & OKE & OXY & SLB & VLO & WMB & XOM\\
Financials    & AFG & AFL & AIG & AJG & AON & AXP & BAC & BEN \\
              & BK & BXS &  C & GL & JPM  & L  & LNC & MMC \\
              & PGR & PNC & RJF & SCHW & STT & TROW & TRV & UNM\\
              & USB & WFC & WRB & WTM\\
Health  Care       &   ABMD & ABT & AMGN & BAX & BDX & BIO & BMY & CAH\\
          & CI & COO & CVS & DHR & HUM & JNJ & LLY & MDT \\
              & MRK & OMI & PFE & PKI & SYK & TFX & TMO & VTRS\\
              & WST\\
Industrials   &   ABM & AIR & ALK & AME & AOS & BA & CAT & CMI\\
              &   CSL & CSX & DE  & DOV & EFX & EMR & ETN & FDX\\
              &   GD & GE & GWW & HON & IEX & ITW & JCI & KSU \\
              & LMT & LUV & MAS & MMM & NOC & NPK&   NSC & PCA\\
              & RPH & PNR & ROK & ROL & RTX & SNA & SWK & TXT\\
              &   UNP\\
Information &   AAPL & ADBE & ADI & ADP & ADSK & AMAT & AMD & GLW \\
Technology  & HPQ & IBM & INTC & MSFT & MSI & MU & ORCL & ROG \\
            & SWKS & TER & TXN & TYL& WDC & XRX\\
Materials  &   APD & AVY & BLL & CCK & CRS & ECL & FMC & GLT \\
          & IFF & IP &  MOS & NEM & NUE & OLN & PPG & SEE  \\
          & SHW & SON & VMC\\
Real Estate      &   ALX & FRT & GTY & HST & PEAK & PSA & VNO & WRI\\
    & WY \\
Utilities   & AEP &  ATO  & BKH & CMS & CNP & D & DTE & DUK\\
            & ED & EIX & ETR & EVRG & EXC & LNT & NEE & NFG \\
            & NI & NJR & OGE &  PEG & PNM & PNW & PPL & SJW\\
            & SO & SWX & UGI & WEC & XEL\\
\hline
    \end{tabular}
\caption{U.S. stock data, sorted by GICS Sector. Data provided by Quandl.}
\label{table:stocks}
\end{table}

\begin{table}[h]
\centering
\begin{tabular}{ l l l l l}
\hline
\textbf{Asset Class} &  & & \textbf{Market (Symbol)} &  \\
\hline
Energy & & WTI crude (CL)       & Heating oil (HO) & Gasoil (QS)\\
      &  & RBOB gasoline (XB)   &                  &            \\
\\
Grain & & Bean oil (BO)         & Corn (C)         & KC Wheat (KW)\\
      & & Soybean (S)           & Soy meal (SM)    & Wheat (W)\\
\\
Livestock & & Feeder cattle (FC) & Live cattle (LC) &  Lean hogs (LH)\\
\\
Metal   &  & Gold (GC)          & Copper (HG)      & Palladium (PA)\\
        &  & Platinum (PL)      & Silver (SI)      & \\
\\

Soft   &   & Cocoa (CC)         & Cotton (CT)      & Robusta Coffee (DF)\\
       &  & Coffee (KC)        & Canola (RS)       & Sugar (SB)  \\

\hline
\end{tabular}
\caption{Futures market data, following Bloomberg symbology. Data provided by Commodity Systems Inc.}
\label{table:futures}
\end{table}

\end{document}